\documentclass{article}
\usepackage{graphicx,amsmath,amssymb,amsfonts,color,fullpage} 

\numberwithin{equation}{section}

\newcommand{\id}{{\rm id}}
\newcommand{\Cc}{{\mathbb C}}
\newcommand{\Oc}{\mathcal O}
\newcommand{\GL}{\mathrm{GL}}

\newcommand\SFT{\mathcal S}
\title{The versatility of the Drinfeld double of a finite group}
\author{Giovanna Carnovale, Nicola Ciccoli, Elena Collacciani}
\date{}

\begin{document}

\maketitle
\hskip5cm\emph{Voor Tom, voor zijn speciale functie binnen de wiskundige gemeenschap}

\begin{abstract}
The Drinfeld double of a finite groups appears in many different areas of mathematics and physics. 
    We review different instances in which the Drinfeld double of a finite group and its representations play a role, focusing on some of Tom Koornwinder's research interests: harmonic analysis, Lie algebras, quantum groups, non-commutative geometry, and Verlinde formula for fusion rules. 
\end{abstract}
\section{Introduction}

In order to show that the category of finite-dimensional representations of  a quantized enveloping algebra $U_q$ is braided monoidal, Drinfeld showed that $U_q$ is a quotient of a
 (topological) quasitriangular Hopf algebra $D({U}^{\geq0})$, called the {\em quantum double}, constructed from a quantized Borel subalgebra $U^{\geq0}$, \cite[Section 13]{drinfeld}.
 As underlying vector space, $D({U}^{\geq0})$ is the (topological) tensor product of $U^{\geq0}$ and an algebra  $U^\circ$ in duality with it, and contains $U^{\geq0}$ and the co-opposite of $U^\circ$ as Hopf subalgebras. It was the first appearance of what is now called the {\em Drinfeld double} of a Hopf algebra. The construction was inspired by an analogous one naturally appearing in the context of Poisson-Lie groups. 
 
\bigskip

 Containing information about both a Hopf algebra $H$ and its dual $H^*$, the Drinfeld double $D(H)$ of a Hopf algebra $H$ is very close to $D(H^*)$, and its category of representations is rich of symmetry. This is also visible through the many subcategories of modules of $H$ that are equivalent to the category of modules of $D(H)$.

\bigskip
 
 Even the simplest example, namely the Drinfeld double  $D(G)$ of the group algebra of a finite group $G$ is already extremely rich as it has connections with many different areas of mathematics and physics, ranging from algebra to analysis and from geometry to gauge theory. Already in the early 1990es, it became apparent that the Drinfeld double of relevant symmetry groups, its representations, and different types of Fourier transforms related to it were underlying many different constructions and theories ranging from representations of finite groups of Lie type in the work of Lusztig \cite{Orange},  link and knot invariants and representations of the braid group \cite{gould}, the classification of Hopf algebra structures on path algebras \cite{CR}, quantum algebras of observables \cite[Chapter VI]{mabook}, conformal field theories and the Verlinde formula \cite{DPP}. It has been frequently used in the theoretical physics literature as mentioned in \cite{muger}, whose primary goal was to provide a mathematically rigorous construction of quantum field theories with $D(G)$-symmetry for any finite $G$. 
 
\bigskip

The double $D(G)$ is extremely versatile: it can be seen as an algebra of functions on $G\times G$ with extra symmetry; its representations can be viewed as $G$-equivariant vector bundles or as $G$-equivariant sheaves on the variety $G$ with finitely many points. The gained analytic or geometric insight led to the above mentioned discoveries.

\medskip

For instance, the ribbon Hopf algebra structure of $D(G)$ is used to construct representations of  mapping class groups of surfaces, \cite{ly,tu}. The case of the torus is especially relevant, as it allows to define a $\mathrm{SL}_2(\mathbb Z)$-action on the span of characters of $D(G)$. The action of the generator $s$ in the usual presentation with relations $s^4=1$ and $s^2=(st)^3$ is  interpreted as a Fourier transform \cite{LM,tom, exotic, CW} and it allows to give an algebraic proof of Verlinde fusion rules, \cite{tom,V}, see also the detailed account in \cite{Br}. A similar Fourier transform for $D(G)$ was defined through the matrix encoding a change of basis between irreducible characters and almost-characters in Lusztig's program for the description of the character tables of finite groups of Lie type in a uniform way, \cite{Orange,arcata,exotic}. 

\bigskip

Ever since,  $D(G)$, its twistings and their representations were used in a broad variety of situations including among others the production of models for quantum computing \cite[Section 5]{kitaev},\cite[Section 6]{rowell},\cite[Section 2.1]{rowell2} and Andruskiewitsch and Schneider's program of classification of finite-dimensional pointed Hopf algebras \cite{AS}. Even at the date of today, almost 40 years after its introduction in Drinfeld's notes from the ICM1986, new interpretations and applications to physics and to novel theories appear, for instance in Majid and McCormak's recent preprint \cite{majmc}. It is shown there how its irreducible representations arise as a non commutative analogue of Mackey quantization of a space with a finite number of points and non trivial internal symmetries and it is hinted that this could play a role in understanding Kitaev's work on quantum computing. 
The importance of the Drinfeld double of a group is nowadays accepted as to deserve a full chapter in basic algebra textbooks \cite[Chapter 18]{Br}.

\bigskip

Just as groups and their representations encode all different types of symmetry,  a pervasive concept in mathematics, physics, and natural sciences, it is natural that also the representations of $D(G)$, encoding extra symmetry, occur in a huge amount of contexts and we cannot hope to provide an exhaustive account of all existing applications and properties of $D(G)$ and its modules in a single paper. In this small survey we give a glimpse of some of the different points of view on $D(G)$ from various areas in mathematics, touching upon some of Tom Koornwinder's research interests: harmonic analysis, Lie algebras, quantum groups, non-commutative geometry, and Verlinde formula.

\section{Notation}
Unless stated otherwise, throughout this survey $k$ is a field and unadorned tensor products are to be considered over $k$. Unless otherwise stated Hopf algebras will be assumed to be finite-dimensional hence they have bijective antipode. For any Hopf algebra $\mathcal H$, we denote by $\mathcal H^{op}$ its opposite Hopf algebra and by $\mathcal H^{cop}$ its co-opposite Hopf algebra.  We adopt Sweedler-like notation with implicit summation  that is, for $\Delta$ a coproduct on $\mathcal H$ and $h\in \mathcal H$ we write $\Delta(h)=h_{(1)}\otimes h_{(2)}$ and for $V$ a right $\mathcal H$-comodule with comodule map $\rho\colon V\to V\otimes \mathcal H$ and $v\in V$ we write $\rho(v)=v_{(0)}\otimes v_{(1)}$ and similarly for their appropriate iterations.   The category of 
left (respectively, right) $\mathcal H$-modules is denoted by $_{\mathcal H}\mathcal M$ (respectively $\mathcal M_{\mathcal H}$) and the category of left (respectively, right) comodules will be denoted by $^{\mathcal H}\mathcal M$ (respectively $\mathcal M^{\mathcal H}$). If $k=\mathbb C$ and $\mathcal H$ is a Hopf $\star$-algebra, we denote by $h^\star$ the image of $h\in\mathcal H$ through the $\star$-operation.
For any $n$-tuple of vector spaces $V_1,\,\ldots,\,V_n$, and any pair of indices $i,\,j\in\{1,\,\ldots,\,n\}$ with $i<j$ we denote by $\tau_{ij}$ the isomorphism $V_1\otimes\cdots\otimes V_i\otimes \cdots\otimes V_j\otimes \cdots\otimes V_n\to V_1\otimes\cdots\otimes V_j\otimes \cdots\otimes V_i\otimes \cdots\otimes V_n$ interchanging the $i$-th and $j$-th tensorands. When $n=2$ we simply write $\tau$. If $f\colon V\to W$ is a linear map of vector spaces, we denote by $f^*$ the adjoint (dual) map $W^*\to V^*$. 

For expository reasons, we will omit the associator in monoidal categories, which in the treated cases comes from natural identifications of vector spaces. 

Unless otherwise stated  $G$ indicates a finite group, $kG$ its group algebra, viewed as a Hopf algebra in the standard way and $C(G)$ the Hopf algebra dual to $kG$. For $g\in G$, we denote the centraliser of $g$ in $G$ by $C_G(g)$. 
If $G$ acts on a set $X$, we denote  the orbit set by $X/G$. 

\section{The Drinfeld double}\label{sec:basics}

 Let $H$ be a  Hopf algebra over a field $k$, with product $m_H$, coproduct $\Delta$, unit element $1_H$, counit $\varepsilon_H$, and antipode $S_H$, and let $H^*$ be its dual Hopf algebra, with operations $\Delta^*$, $m_H^*$, unit $\varepsilon_H^*$, counit $(-,1_H)$, and antipode $S_H^*$. The Drinfeld double $D(H)$ of $H$ is uniquely determined by the following conditions:

 \begin{enumerate}
 \item[I.] As a vector space, $D(H)=H^*\otimes H$;
\item[II.] $D(H)$ is a Hopf algebra;
     \item[III.]\label{item:Hopf-sub} $H\simeq k\varepsilon_H\otimes H$ and $(H^*)^{cop}\simeq (H^*)^{cop}\otimes k1_H$ are Hopf subalgebras of $D(H)$;
     \item[IV.] For any $f\in H^*$ and any $h\in H$\begin{equation}\label{item:compa}(\varepsilon_H\otimes h)(f\otimes 1_H)=f(S_H^{-1}(h_{(3)})-h_{(1)})\otimes h_{(2)}.\end{equation} 
     \end{enumerate}
The antipode $S_{D(H)}$ is given by 
$S_{D(H)}(f\otimes h)=f\left(h_{(3)} S_H^{-1}(-)S_H^{-1}(h_{(1)})\right)\otimes S_H(h_{(2)})$ for $f\in H^*$ and $h\in H$ and
it is again bijective by construction. In Hopf theoretical terms $D(H)$ can be seen as a cocycle twist of $H^{*,cop}\otimes H$ or as a double-crossed product for the mutual coadjoint actions of $H$ and $H^*$, \cite[Proposition 7.1.5, Example 7.2.5]{mabook}. We recall that the {\em coadjoint action} of a Hopf algebra $H$ on its dual $H^*$ is given by 
\begin{equation}\label{eq:coadjoint}(Ad^*(h)\cdot f)=f(S_{H}(h_{(1)})\,-\,h_{(2)})\mbox{ for }h\in H,\mbox{ and }f\in H^*.\end{equation}

A key property of  $D(H)$ is that it is quasi-triangular \cite{drinfeld} with universal $R$-matrix $R_{D(H)}$ given by the canonical element, that is, for any dual pair of bases $\{e_i\}_{i\in I}$ and  $\{e^i\}_{i\in I}$ of $H$ and $H^*$, respectively,  \begin{equation*}
R_{D(H)}=\sum_{i\in I}(\varepsilon\otimes e_i)\otimes (e^i\otimes 1_H).
\end{equation*} 
Recall that for any $R$-matrix $\mathcal R$ in a quasitriangular Hopf algebra $\mathcal H$, the element  $(\tau \mathcal R)^{-1}\in \mathcal H\otimes \mathcal H$ is also an $R$-matrix for $\mathcal H$, \cite[Lemma VIII.2.3]{kassel}, so  
\begin{equation*}
R'_{D(H)}=(\tau R_{D(H)})^{-1}=\sum_{i\in I}(e^i\otimes 1_H)\otimes(\varepsilon\otimes S_H(e_i))
\end{equation*}
is an $R$-matrix for $D(H)$. In fact, $D(H)$, with any of the above $R$-matrices is {\em minimal quasitriangular}, that is, it contains no quasitriangular Hopf subalgebras. Any minimal quasitriangular Hopf algebra is the quotient of a Drinfeld double, \cite{radford}.   

Through the literature, the Drinfeld double is realized also on the underlying vector space $H\otimes H^*$, and/or having 
 $(H^*)^{op}$ as a Hopf subalgebra. The four different realizations are given in  \cite[Thereom 7.1.1, Exercise 7.1.2]{mabook}. We have adopted here the formulation from \cite{kassel,radbook}, for uniformity with other used results.

\medskip

To fix notation, we recall that if a Hopf algebra $\mathcal H$ is quasitriangular with universal $R$-matrix $\mathcal R$, then the category $_{\mathcal H}\mathcal M$ of its left modules is braided: for any pair of  $\mathcal H$-modules $V$ and $W$, there is a natural $\mathcal H$-module isomorphism 
$c_{V,W}\colon V\otimes W\to W\otimes V$, called {\em braiding}, defined on pure tensors $v\otimes w$  by $c_{V,W}(v\otimes w)=\tau \circ (\mathcal R(v\otimes w))$. The properties of the $R$-matrix ensure that for each module $V$, the isomorphism $c_{V,V}$ satisfies the Yang-Baxter equation
\begin{equation}\label{eq:YB}
(c_{V,V}\otimes \id_V)(
\id_V\otimes c_{V,V})(c_{V,V}\otimes \id_V)=(
\id_V\otimes c_{V,V})(c_{V,V}\otimes \id_V)(
\id_V\otimes c_{V,V}).
\end{equation}
The quest for solutions of this equation was one of the motivations of quantum group theory. 

\bigskip

We will also need two special elements associated to a quasitriangular Hopf algebra $(\mathcal H,\mathcal R)$: the product  $Q:=(\tau \mathcal R)\mathcal R\in \mathcal H\otimes\mathcal H$,  called {\em monodromy element} and $u:=m(S_{\mathcal H}\otimes\id)(\tau \mathcal R)\in\mathcal H$, called the {\em Drinfeld element} of $\mathcal H$. In the terminology introduced in \cite{reshe}, the quasitriangular Hopf algebra $(D(H),R_{D(H)})$ is also {\em factorizable}, that is, the linear map 
\begin{equation}\label{eq:fQ}
f_Q\colon \mathcal D(H)^*\to\mathcal D(H),\quad \zeta\mapsto (\zeta\otimes\id_{\mathcal H})Q\end{equation} is an isomorphism.

\medskip

By construction, it is natural to expect the double to have plenty of symmetries. A result in this sense is \cite[Theorem 3]{radford}.
It states that, using the identification $H^{**}\simeq H$, we have a Hopf algebra isomorphism
\begin{equation}\label{eq:dual}T:=\tau(S_{H^*}^{-1}\otimes S_H)\colon D(H)\longrightarrow D(H^*)^{op} \quad\mbox{ satisfying }\quad(T\otimes T)(R_{D(H)})=\tau R_{D(H^*)}.\end{equation}

We remind  here that by standard theory, $D(H)^{op}$ and $D(H)^{cop}$ are also quasi-triangular, with universal $R$-matrix $R_{D(H)}^{-1}$.

\bigskip

Assume now that $k={\mathbb C}$. Then one gains the possibility of using the machinery of $\star$-operations on Hopf algebras, allowing for the definition of real forms. The idea of introducing a $\star$-operation on an algebra with what we now call a coproduct is probably firstly due to G. I. Kac, who developed the notion of ring group in \cite{Kac}. This should be seen as a Hopf algebra in the context of von Neumann algebras, see \cite{AFS}. In order to develop a theory of compact quantum groups, a star operation in Hopf algebras in the context of $C^*$-algebras was introduced in \cite{30}. According to \cite[Section 2]{maj-qm} the $\star$-structure is also useful for addressing problems in quantum mechanics, as it allows to consider representations on Hilbert spaces, and evaluate through $\star$ the positivity of states of observables. 

\medskip

Formally, a Hopf $\star$-algebra is the datum of: a Hopf algebra $H$ together  with a $\star$-operation, that is,
 a complex conjugate-linear, anti-multiplicative involution $\star \colon H\to H$ denoted on elements by  $h\mapsto h^\star$, preserving the unit element, satisfying $\Delta\circ \star=(\star\otimes \star)\circ\Delta$, and such that $\varepsilon$ intertwines $\star$ and complex conjugation. These conditions imply that $(S_H\circ \star)^ 2=\id$. If $H$ is a Hopf $\star$-algebra, then $H^*$ is also a Hopf $\star$-algebra, with $\star$-structure given by $f^\star(h)=\overline{f(S_H(h)^\star)}$ for $h\in H$ and $f$ in $H^*$. By \cite[Example 2.1]{maj-qm}  if $H$ is a Hopf $\star$-algebra, then $D(H)$ is also a Hopf $\star$-algebra and  $R_{D(H)}^{\star\otimes \star}=R_{D(H)}^{-1}$, that is, it is anti-real.
 We point out that in principle one may define a $\star$-structure on Hopf algebras over any field equipped with an involution playing the role of complex conjugation. 

\bigskip 

Our main focus is the case in which $H=kG$ is the group algebra of a finite group $G$, we write $D(G)$ instead of $D(kG)$ for simplicity. It can be described in different ways, according to which of its features one wishes to underline. We collect here some of its incarnations. 

\subsection{$D(G)$ for the algebraist}\label{sec:DG-alg}
A basis for $kG^*$ is given by the delta functions $\delta_g$ for $g\in G$, that is, $\delta_g(g')=\delta_{g,g'}$ for any $g'\in G$.
For $g,h,l\in G$, we have 
\begin{equation}\label{eq:prod-H^*}
(\delta_g\delta_h)(l)=\Delta^*(\delta_g\otimes\delta_h)(l)=\delta_{g,l}\delta_{h,l}=\delta_{g,h}\delta_g(l),
\end{equation}
that is, the elements $\delta_g$ for $g\in G$ are a set of central idempotents for $(kG)^*$, and $1_{kG^*}=\varepsilon_{kG}=\sum_{g\in G}\delta_g$. 
In addition, given $g,h,l\in G$, we have
\begin{equation*}
\Delta_{kG^*}(\delta_g)(h\otimes l)=m_H^*\delta_g(h\otimes l)=\delta_g(hl)=\delta_{g,hl}=\sum_{t,s\in G, ts=g}(\delta_s\otimes\delta_t)(h\otimes l)
\end{equation*}
so $\Delta_{kG^*}(\delta_g)=\sum_{t,s\in G, ts=g}(\delta_s\otimes\delta_t)$ and for $g,h\in G$ we have
\begin{equation}
S_{kG^*}(\delta_g)(h)=\delta_g(h^{-1})=\delta_{g,h^{-1}}=\delta_{g^{-1}}(h).
\end{equation}

Omitting the $\otimes$ symbol for simplicity, a basis for $D(G)$ is given by $\delta_hg$, for $h,\,g\in G$. Combining \eqref{eq:prod-H^*} with condition IV gives the multiplication in $D(G)$:
\begin{equation}
(\delta_hg)(\delta_t l)=\delta_h\delta_{gtg^{-1}}gl=\delta_{g^{-1}hg,t}\delta_h gl
,\quad \textrm{ for }g,h,l,t\in G
\end{equation}
and $1_{D(G)}=\sum_{g\in G}1_G\delta_g$. Condition III gives the comultiplication:
\begin{equation}
\Delta_{D(G)}(\delta_hg)=\sum_{s,t\in G, st=h}\delta_tg\otimes \delta_sg,\quad \mbox{ for }g,h\in G. 
\end{equation}
Since the antipode is an algebra antimorphism and $S_{(kG^*)^{op}}=S_{kG^*}^{-1}$, there holds
\begin{equation}
S_{D(G)}(\delta_hg)=S_{kG}(g)S_{kG^*}^{-1}(\delta_h)=g^{-1}\delta_{h^{-1}}=\delta_{g^{-1}h^{-1}g}g^{-1}
   \textrm{ for }g,h\in G.\end{equation}
Finally, counit and $R$-matrices are
\begin{align}\label{eq:errematrice}
\nonumber&\varepsilon_{D(G)}(\delta_h g)=\varepsilon_{kG^*}(\delta_h)\varepsilon_{kG}(g)=\delta_h(1_G)=\delta_{h,1}, \textrm{ for }g,\,h\in G \\
&R_{D(G)}=\sum_{h,g\in G}\delta_h g\otimes \delta_g1_G \textrm{ and } R'_{D(G)}=\sum_{h,g\in G}\delta_g1_G \otimes \delta_hg.
\end{align}

\bigskip

In general $R_{D(H)}$ and $R'_{D(H)}$ may not be the unique $R$-matrices on $D(H)$: for $G$ a finite group the quasitriangular structures on $\mathcal H= D(G)$ have been classified in \cite{Ke}. 

\bigskip

Several authors have computed the center of $D(G)$,  see for instance \cite[Section 8.1.3]{br}. It can be proved  by direct calculation to be the subalgebra spanned by the orbit sums 
\begin{equation}\label{eq:center}
\Sigma_{G(h,g)}:=\sum_{l\in G}\delta_{lhl^{-1}}lgl^{-1},\quad\mbox{ for } h,\,g\in G,\;[h,g]=1.
\end{equation}  

The Hopf algebra structure of $D(G)$ somehow determines the group $G$, in the sense that for each pair of groups $G_1,\,G_2$ the isomorphism $D(G_1)\simeq D(G_2)$ implies $G_1\simeq G_2$, \cite{ke2}.
 
\bigskip

If $k=\Cc$, the natural $\star$-structures on $\Cc G$ and $\Cc G^*$ respectively, are given by the conjugate-linear extensions of the assignments $g^\star=g^{-1}$ and $\delta_g^\star=\delta_g$ for any $g\in G$. They extend to a $\star$-operation on $D(G)$ by setting
\begin{equation}\label{eq:star-algebra}
(\delta_hg)^\star=\delta_{g^{-1}hg} g^{-1},\quad \forall g,h\in G.
\end{equation}

\subsection{D(G) for the analyst}\label{sec:analyst}

In \cite{tom} the authors describe $D(G)$  through an analyst's lens: their apparently different description is best suitable for the generalization to the case of compact and locally compact Lie groups carried out in \cite{tom-tensor,tom-muller}. We recall here their approach. Observe that the assumption $k=\mathbb C$ can be relaxed, as the assumption that ${\rm char} (k)$ does not divide $|G|$ suffices.  Notice that in \emph{loc. cit.}  the symbol $D(G)$ stands for $D(kG^*)$, but the underlying vector space is as in \cite{CP}, so it matches with the one we used for $D(G)$. As a matter of fact, in virtue of the isomorphism \eqref{eq:dual}, the Hopf algebra $D(kG^*)$ is isomorphic to $D(G)^{op}$, which, through $S_{D(G)}$ is isomorphic to $D(G)^{cop}$. Hence, it is equivalent to study the algebra representations of $D(kG^*)$ or of $D(G)$. 

\bigskip

The finite group $G$ is naturally endowed with the discrete topology and has a normalized invariant measure:
\[\int_Gf(g)dg=\frac{1}{|G|}\sum_{g\in G}f(g).\]
For uniformity with Section \ref{sec:DG-alg}, the delta function has weight $1$, differently than in \cite{tom}, that is, we use a Kronecker delta rather than a Dirac's delta. 

As a vector space $kG$ can be identified with the space $C(G)$ of (continuous) functions on $G$ by mapping the linear combination $\sum_g a_g g$ to the function $a\colon G\to k$ mapping each $g\in G$
to its coefficient $a_g$. In other words, we map $g$ to $\delta_g$. Through this identification, the product of the functions $a$ and $b$ is the function $(a\cdot b)$ whose evaluation at $g\in G$ is $\sum_{h\in G}a(h)b(h^{-1}g)$. The neutral element is then the delta function $\delta_{1_G}$ and the counit is $f\mapsto |G|\int_G fdg$. 

We can also identify the vector space $kG^*$ with $C(G)$. Through this identification, the product in $kG^*$ is pointwise multiplication of functions, the coproduct  $\Delta_{C(G)}\colon C(G)\to C(G)\otimes C(G)\simeq C(G\times G)$ is given by $\Delta_{C(G)}(f)(g,h)=f(gh)$, for $f\in C(G)$ and $g,h\in G$ and the counit $\varepsilon_{C(G)}$ is evaluation at $1_G$. 

\medskip

For compatibility between the presentation in Subsection \ref{sec:DG-alg} and the one in \cite{tom} we choose the chain of vector space identifications  
\begin{equation}\label{eq:identify}
D(G)=kG^*\otimes kG=C(G)\otimes kG\simeq C(G)\otimes C(G)\simeq 
 C(G\times G)\end{equation} 
 where the basis element $\delta_hg\in kG^*\otimes kG$ corresponds to the delta function $\delta_{(h^{-1},g)}$ on $G\times G$. Similarly,  we  have the chain of identifications
\begin{equation} D(G)\otimes D(G)\simeq C(G\times G)\otimes C(G\times G)\simeq C((G\times G)\times (G\times G))\end{equation} whose composition maps 
$\delta_{h_1}g_1\otimes\delta_{h_2}g_2$ to $\delta_{(h_1^{-1},g_1,h_2^{-1},g_2)}$. Then multiplication, unit, counit, antipode and universal $R$-matrices 
from Section \ref{sec:DG-alg} become:
\begin{align*}
&(f_1\cdot f_2)(h,g)=|G|\int_G f_1(h, g_1)f_2(g_1^{-1}hg_1,g_1^{-1}g)dg_1,&&1_{D(G)}(h,g)=\delta_{1_G}(h),\\
&(\Delta(f))(h_1,g_1,h_2,g_2)=f(h_1h_2,g_1)\delta_{g_1, g_2},&&\varepsilon(f)=|G|\int_Gf(1_g,g)dg,\\
&S(f)(h,g)=f(g^{-1}h^{-1}g,g^{-1}),\\
&R(h_1,g_1,h_2,g_2)=\delta_{1_G}(g_1h_2)\delta_{1_G}(g_2),
&&R'(h_1,g_1,h_2,g_2)=\delta_{1_G}(g_1)\delta_{1_G}(h_1g_2^{-1})
\end{align*}
for $h,g,h_1,h_2,g_1,g_2\in G$.

\medskip

If $k=\Cc$, then the $\star$-structure \eqref{eq:star-algebra} through the identification \eqref{eq:identify} becomes $f^\star(h,g)=\overline{f(g^{-1}hg,h^{-1})}$
for $f\in C(G\times G)$ and $g,h\in G$. In addition, $D(G)$ and $D(G)^*$ have faithful positive linear Haar functionals, and they are $C^*$-algebras \cite[Section 2]{tom}, so they both fit in the framework of compact matrix pseudogroups of \cite{W}. 

\bigskip

Other realisations of $D(G)$ as an algebra of functions have been considered in the literature: for instance, \cite[(9)]{tom-muller} uses contraction to identify the underlying vector space of $D(G)$ with $C(G, kG)$, by sending $f\otimes g\in C(G)\otimes kG$ to the function $z\mapsto f(z)g$.

\section{The different faces of the category of $D(G)$-modules}

The category of left $D(H)$-modules is braided monoidal, with braiding given by the action of $\tau R_{D(H)}$. It encodes much of the Hopf algebra structure of $H$. In the early 90's several equivalences of (braided monoidal) categories of modules have been discovered, establishing connections between a variety of Hopf algebraic objects that had been introduced independently.  We first spell out the explicit case of $H=kG$ and then recall some of these connections.

\medskip 

Let $V$ be a left $D(H)$-module with action denoted by the $\cdot$ symbol. Then $V$ is a left $H$-module and a left $H^*$-module by restriction, and the two actions must statisfy the constraint stemming from  \eqref{item:compa}:
\begin{equation}\label{eq:D-action}
    ((\varepsilon\otimes h)(f\otimes 1_H))\cdot v=(f(S_H^{-1}(h_{(3)})-h_{(1)})\otimes h_{(2)})\cdot v, \quad\mbox{ for } h\in H,\,f\in H^*,\, v\in V.
\end{equation}

\medskip

When $H=kG$, the condition of $V$ being a left $kG^*$-module is equivalent to the requirement that $V$ is $G$-graded, that is, $V=\bigoplus_{g\in G}V_g$ where $V_h=\delta_h(V)$ and the idempotent $\delta_h$ for $h\in G$ acts as the projection onto $V_h$ along the given decomposition of $V$. Then, \eqref{eq:D-action} becomes 
\begin{equation}\label{eq:DG-mod}gV_h=V_{ghg^{-1}},\quad \mbox{ for }g,h\in G.\end{equation}

\medskip

Thus, the objects in the category $_{D(G)}{\mathcal M}$ are $G$-graded $kG$-modules satisfying \eqref{eq:DG-mod}, and the morphisms are the $G$-equivariant $G$-graded maps. The tensor product of objects $V=\bigoplus_{g\in G}V_g$ and $W=\bigoplus_{g\in G}W_g$ in  $_{D(G)}{\mathcal M}$ is $V\otimes W$ with usual $G$-action and with $G$-grading satisfying $(V\otimes W)_h=\sum_{g\in G}V_{g^{-1}h}\otimes W_g$ for $h\in G$. The braiding  $c_{V,W}\colon V\otimes W\to W\otimes V$ is given by 
\begin{equation} \label{eq:braiding}
v\otimes w\mapsto w\otimes h\cdot v,\quad \forall v\in V,\;h\in G,\; w\in W_h. 
\end{equation} 

The {\em support} ${\rm supp}(V)$ of the $G$-grading in a $D(G)$-module $V$ is an invariant of $V$ that is stable by conjugation by $G$. In addition, $D(G)$-modules with disjoint supports have only trivial morphisms between them. Therefore, denoting by $_{D(G)}{\mathcal M}^{\Oc}$ the subcategory of 
$_{D(G)}{\mathcal M}$ whose objects are supported in the conjugacy class $\Oc$ of $G$, we have a 
block decomposition of categories
\begin{equation}\label{eq:deco-category}
_{D(G)}{\mathcal M}=\bigoplus_{\Oc_g\in G/G} {}_{D(G)}{\mathcal M}^{\Oc_g}
\end{equation}
see also \cite[Proposition 3.7]{CR},\cite[Proposition 18.1.15]{Br}. For two conjugacy classes  $\Oc_1,\,\Oc_2\in G/G$ the tensor product satisfies
\begin{equation*}
\otimes\colon\; _{D(G)}{\mathcal M}^{\Oc_2} \times {}_{D(G)}{\mathcal M}^{\Oc_1} \to \bigoplus_{\Oc\in G/G, \Oc\cap(\Oc_1\Oc_2)\neq\emptyset}\;{}_{D(G)}{\mathcal M}^{\Oc}\end{equation*}
and the braiding $c$ satisfies
\begin{equation*}c(_{D(G)}{\mathcal M}^{\Oc_1} \otimes{}_{D(G)}{\mathcal M}^{\Oc_2})\subseteq\; _{D(G)}{\mathcal M}^{\Oc_2} \otimes\;{}_{D(G)}{\mathcal M}^{\Oc_1}\end{equation*}
in virtue of \eqref{eq:braiding} and \eqref{eq:DG-mod}.

\subsection{Yetter-Drinfeld modules}
Let $V$ be a left $D(H)$-module. The left $H^*$-module structure is equivalent to a right $H$-comodule structure on $V$, where $f\cdot v=v_{(0)}f(v_{(1)})$ for $v\in V$ and $f\in H^*$. It was observed in \cite[Prop 2.2]{Maj} that \eqref{eq:D-action} is equivalent to the condition 
\begin{equation}\label{eq:YD-module}
h_{(1)}\cdot v_{(0)}\otimes h_{(2)}v_{(1)}=(h_{(2)}\cdot v)_{(0)}\otimes (h_{(2)}\cdot v)_{(1)}h_{(1)},
\end{equation}
see also \cite[Theorem IX.5.2]{kassel} for a detailed proof. 

\medskip

A left $H$-module and right $H$-comodule satisfying \eqref{eq:YD-module} is nowadays called a {\em (left-right) Yetter-Drinfeld module}, since a slightly different but essentially equivalent version of this condition firstly occurred in the work of Yetter \cite[Definition 3.6]{yetter}, where such structures were called {\em crossed bimodules}. More precisely, there are $4$ versions of Yetter-Drinfeld modules, according to whether the module and comodule structures of $H$ are left or right: denoting respectively by 
$^H\mathcal{YD}_H$,  $^H_H\mathcal{YD}$, $\mathcal{YD}^H_H$ and $_H\mathcal{YD}^H$ the categories  whose objects are respectively right-left, left-left, right-right, and left-right Yetter-Drinfeld modules and whose morphisms are the $H$-linear and $H$-colinear maps, we have the following equivalences of categories, \cite[Proposition 6]{rad-tow}:
\begin{equation}\label{eq:first-equivalences}
^H\!\mathcal{YD}_H\simeq\,  ^H_H\!\mathcal{YD}\simeq\mathcal{YD}^H_H\simeq\,_H\!\mathcal{YD}^H.\end{equation}
Each of these categories is equivalent to $_{D(H)}{\mathcal M}$ in virtue of \cite[Prop 2.2]{Maj}. The four variants essentially correspond to the four different equivalent realisations of $D(H)$ and Yetter's original crossed bimodules are left $H$-modules and left $H$-comodules. When $H$ is cocommutative, as it is for $kG$, its left and right comodules coincide so the number of different variants is reduced. 

Yetter-Drinfeld modules are also called {\em Yang-Baxter modules} \cite{lam-rad}, {\em Yetter-Drinfed structures} in \cite{rad-tow} and {\em crossed $H$-modules} in \cite{mabook}.

The categories of Yetter-Drinfeld modules can be naturally equipped with a braided monoidal structure \cite[Theorem 7.2]{yetter} and the equivalences above extend to braided monoidal equivalences, \cite[Section 2]{Maj}. This way the braided monoidal category $(_{D(H)}{\mathcal M},\otimes,R)$ can be described in terms of $H$ only. 

\medskip

Back to our focus case $H=kG$,  an object  in $^{kG}{\mathcal YD}_{kG}$ is a $G$-module $V$ with a right $kG$-comodule structure $\rho\colon V\to V\otimes kG$, satisfying  
\begin{equation}\label{eq:YD-group}
h\cdot v_{(0)}\otimes hv_{(1)}=(h\cdot v)_{(0)}\otimes (h\cdot v)_{(1)}h
\end{equation} for $h\in G$ and $v\in V$. Now, having $\rho$ is equivalent to giving
a $G$-grading $V=\bigoplus_{g\in G}V_g$, with $\rho(v)=v\otimes g$ for any $v\in V_g$. Moreover, \eqref{eq:YD-group} applied to $v\in V_g$ gives
\begin{equation}\label{eq:yd-group}
h\cdot v\otimes hgh^{-1}=(h\cdot v)_{(0)}\otimes (h\cdot  v)_{(1)},
\end{equation}
that is, $h\cdot v$ is homogeneous and lies in  $V_{hgh^{-1}}$, recovering  \eqref{eq:DG-mod} as expected.

\subsection{Hopf modules}\label{sec:bimodules}

Long before the introduction of the quantum double and of Yetter-Drinfeld modules, a different category of  $H$-modules and comodules with a compatibility condition with the two structures had been considered, namely the category of {\em Hopf modules}, \cite[page 82]{sweedler}. 

\medskip

A right-right Hopf module is a right $H$-module and a right $H$-comodule $V$ such that 
the coaction map is a right module morphism from $V$ to the the tensor module $V\otimes H$. 

Similarly, one defines the notion of right-left, left-right, and left-left Hopf module, according to the side on which $H$ acts and coacts, \cite[\S 3]{Sch}. We denote respectively by ${\mathcal M}_H^H$, $^H{\mathcal M}_H$, $_H{\mathcal M}^H$ and $_H^H{\mathcal M}$ the categories whose objects are respectively right-right, right-left, left-right, and left-left Hopf modules and whose morphisms are the $H$-linear and $H$-colinear maps.    

\medskip

One considers then the subcategory of two-sided two-cosided Hopf modules $_H^H{\mathcal M}_H^H$. Its objects are two-sided $H$-modules and two-sided $H$-comodules whose structures ensure that they are simultaneously objects in ${\mathcal M}_H^H$, $^H{\mathcal M}_H$, $_H{\mathcal M}^H$ and $_H^H{\mathcal M}$.
The morphisms in $_H^H{\mathcal M}_H^H$ are the maps that are simulteneously right and left $H$-linear and right and left $H$-colinear. Such modules occurred in the work of Nichols \cite{nichols} and are precisely the {\em bicovariant bimodules} in \cite[Definition 2.1, Definition 2.1]{wor-diff-cal} introduced to study bicovariant differential calculus on quantum groups. 

The category $_H^H{\mathcal M}_H^H$  can be equipped with two distinct but equivalent monoidal category structures \cite[page 87, Corollary 6.1]{Sch} and \cite[Theorem 5.7]{Sch} establishes an equivalence of monoidal categories between $_H^H{\mathcal M}^H_H$ with any of the two monoidal structures and the category ${\mathcal YD}_H^H$ of right-right (and consequently all four versions of) Yetter-Drinfeld modules, with braided monoidal structure stemming from the one in $_{D(H)}{\mathcal M}$ through the equivalences in \cite{Maj} and \eqref{eq:first-equivalences}. The same result was obtained independently in  \cite[Appendix]{AD},\cite{rosso} and was implicitly contained in \cite[Theorem 2.1]{wor-diff-cal}.

In contrast to the equivalence $_{D(H)}{\mathcal M}\simeq\, _H\!{\mathcal YD}^H$, modules corresponding to each other through the equivalence $_H^H{\mathcal M}_H^H\to\, _H\!{\mathcal YD}^H$ do not have isomorphic underlying vector space, since the equivalence is obtained in one direction by taking left coinvariants, and in the other direction by constructing suitable structures on $H\otimes V$ for a Yetter-Drinfeld module $V$.

\medskip

In our focus case, for $V=\bigoplus_{g\in G}V_g\in\; _{D(G)}{\mathcal M}$, one obtains a bicovariant $kG$-module structure on $kG\otimes V$ by taking the right and left $kG$-module structures $\leftharpoonup$ and $\rightharpoonup$, and the right and left $kG$-comodule structures $\xi\colon kG\otimes V\to kG\otimes kG\otimes  V$ and $\rho\colon kG\otimes V\to kG\otimes V\otimes kG$ as follows
\begin{align*}
&(h\otimes v)\leftharpoonup g=hg\otimes g^{-1}v,&g\rightharpoonup(h\otimes v)=gh\otimes v,\\
&\rho(g\otimes v_h)=g\otimes v_h\otimes gh,&\xi(g\otimes v)=g\otimes g\otimes v,&&\mbox{ for }h,g\in G,\,v\in V,\,v_h\in V_h.
\end{align*}
One recovers $V$ by taking the coinvariants $\{w\in kG\otimes V~|~\xi(w)=1\otimes w\}$, with restriction of the $\leftharpoonup$ action and grading from $\rho$, see \cite[A1, A2]{AD}.

\subsection{Doi-Koppinen modules}

Doi and Koppinen \cite{doi,koppinen} independently introduced a generalized variant of Hopf modules that unified different families of modules such as: Hopf modules, relative Hopf modules, and graded modules. Their approach turned out to be useful for the definition of induction of (unitary) representations of (compact) quantum groups, see \cite{KDC} in this volume and references therein. 

\medskip

We recall the definition of this family of modules. We consider triples $(H,A,C)$ where 
 $A$ is a right $H$-comodule algebra with structure map $a\mapsto a_{(0)}\otimes a_{(1)}\in A\otimes H$ for $a\in A$, and $C$ is a left $H$-module coalgebra with action $\bullet$. Then a {\em left-right $(H,A,C)$-Hopf module} is a $k$-vector space $M$ together with a left $A$-module structure with action $\cdot$ and a right $C$-comodule structure $\rho$ with $\rho(m)=m_{[0]}\otimes m_{[1]}$, for $m\in M$ such that 
 \begin{equation}\label{eq:doi-hopf}\rho(a\cdot m)=a_{(0)}\cdot m_{[0]}\otimes a_{(1)}\bullet m_{[1]}\in M\otimes C
\mbox{ for }m\in M, a\in A.\end{equation}
 The category $_A{\mathcal M(H)}^C$ is then the category whose objects are left-right $(H,A,C)$-Hopf modules and whose morphisms are $A$-linear and $C$-colinear morphisms.  Now \cite[Corollary 2.4]{camishe} states that $H$ is naturally: a right $H^{op}\otimes H$-comodule algebra, via $h\mapsto h_{(2)}\otimes S^{-1}(h_{(1)})\otimes h_{(3)}$, and a left $H^{op}\otimes H$-module coalgebra, with action $(h\otimes l)\bullet g=lgh$ for $h,g,l\in H$. Moreover, the categories $_{H}{\mathcal M}(H^{op}\otimes H)^H$ and $_H{\mathcal YD}^H$ are isomorphic. As a consequence, $_{D(H)}{\mathcal M}$ is a Grothendieck category, i.e., it satisfies the conditions in \cite{tohoku}. The equivalence of categories in this case is the identity at the level of  the underlying vector spaces of objects and at the level of morphisms, cf. \cite[Theorem 2.3]{camishe}.

\medskip

For $H=kG$, the right $kG^{op}\otimes kG$-comodule structure of $kG$ is 
given by $g\mapsto g\otimes g^{-1}\otimes g$, and condition \eqref{eq:doi-hopf}
 becomes $\rho(g\cdot m)=g\cdot m_{[0]}\otimes gm_{[1]}g^{-1}$ for $g\in G$ and $m\in M$, so if $\rho(m)=h\otimes m$ for some $h\in G$, i.e., if $m$ is homogeneous of degree $h\in G$, then $g\cdot m$ is homogeneous of degree $ghg^{-1}$, i.e., we recover \eqref{eq:DG-mod} once more.

\subsection{Categorical interpretation of $D(H)$}

 The symmetry of the category of left $D(H)$-modules has a natural interpretation in terms of the center of the category of $H$-modules.  
 Let us first recall the construction of the center of a monoidal category, originally due to Drinfeld (unpublished). It appeared in \cite[Definition 3]{js} and \cite[Theorem 3.3, Example 3.4]{maj-center}.  We follow the exposition in \cite[XIII.4, XIII.5]{kassel}, restricting ourselves to the case of strict monoidal categories for simplicity.

 \medskip

Given a (strict) monoidal category $({\mathcal C},\otimes)$ with unit object $I$, the objects in the {\em Drinfeld center} ${\mathcal Z}({\mathcal C})$ of ${\mathcal C}$ are pairs $(X, c_{-,X})$ where $X$ is an object in ${\mathcal C}$ and $c_{-,X}$ is a natural isomorphism of functors from $(-\otimes X)\colon {\mathcal C}\to{\mathcal C}$ to $(X\otimes -)\colon {\mathcal C}\to{\mathcal C}$ satisfying  
\begin{equation}\label{eq:center-cat}
c_{A\otimes B, X}=(c_{A,X}\otimes\id_B)(\id_A\otimes c_{B,X}),\quad \forall A,B\in{\mathcal C}.
\end{equation}
A morphism in ${\mathcal Z}({\mathcal C})$ from the object $(X, c_{-,X})$ to the object $(Y,c_{-,Y})$ is a morphism $T\colon X\to Y$ in ${\mathcal C}$ such that
$(T\otimes \id_A)c_{A,X}=c_{A,Y}(\id_A\otimes T)$ for any object $A$ in $\mathcal C$. Observe that since $c_{-,X}$ is a natural transformation, 
$(\id_X\otimes T')c_{A,X}=c_{B,X}(T'\otimes\id_X)$ for any morphism $T'\colon A\to B$ in ${\mathcal C}$.

\medskip

The category ${\mathcal Z}(\mathcal C)$ can be equipped with a monoidal category structure: 
indeed, for any pair of objects $(X,c_{-,X})$ and $(Y,c_{-,Y})$ in ${\mathcal Z}(\mathcal C)$, the family of isomorphisms $c_{A,X\otimes Y}:=(\id_X\otimes c_{A,Y})(c_{A, X}\otimes\id_Y)$ where $A$ runs through all objects in $\mathcal C$, combine to give an object $(V\otimes W, c_{-,X\otimes Y})$ in ${\mathcal Z}(\mathcal C)$ and this construction satisfies the pentagon axiom. In addition, the pair $(I,\id)$ gives a unit element satisfying the triangle axiom. In fact, $({\mathcal Z}(\mathcal C),\otimes, (I,\id))$ can also be endowed with a braided monoidal category structure, by setting
$c_{(X, c_{-,X}),(Y,c_{-,Y})}:=c_{X,Y}$. 

\medskip

This construction is universal, in the sense that, calling $\pi\colon \mathcal Z(\mathcal C)\to \mathcal C$ the functor mapping $(X,c_{-,X})$ to $X$ on objects, there exists a unique tensor functor $F\colon {\mathcal T}\to{\mathcal Z}(\mathcal T)$ such that $\pi\circ F=\id_{\mathcal T}$ for any strict braided monoidal category ${\mathcal T}$. 

\medskip

Let now $(V,c_{-,V})$ be an object in $\mathcal Z(_H{\mathcal M})$. Then, $V$ is  an $H$-module by construction. It is also a right $H$-comodule with structure morphism: $\rho(v)=c_{H,V}(1\otimes v)$ for $v\in V$. The two structures satisfy the compatibility condition \eqref{eq:YD-module} and morphisms in $\mathcal Z(_H{\mathcal M})$ are $H$-linear and $H$-colinear. This assignment induces an equivalence of braided monoidal categories \begin{equation}({\mathcal Z}(_H{\mathcal M}),\otimes,  c,k)\to (_H{\mathcal YD}^H,\otimes, c, k)\simeq (_{D(H)}\mathcal{M},\otimes, \tau\circ R_{D(H)}, k)\end{equation} 
such that the restriction monoidal functor $(_{D(H)}\mathcal{M},\otimes, k)\to (_{H}\mathcal{M},\otimes, k)$ stemming from the inclusion $H\subseteq D(H)$ corresponds to $\pi$. 
The connection between $\mathcal Z(_H\mathcal M)$ and the double is attributed to Drinfeld in \cite[Example 3.4]{maj-center}.

\medskip

Dually, $_{D(H)}{\mathcal M}$ is also equivalent to the  Drinfeld center of the category of $H$-comodules $^H{\mathcal M}$, see \cite[Proposition 7.15.3]{egno}. This can  also be seen as a consequence of the equivalence for left $H$-modules and the isomorphism \eqref{eq:dual}, making use of  \cite[Exercises 7.13.5,8.5.2]{egno}, \cite[Lemma VIII.2.3]{kassel}.

\bigskip

 The case of $H=kG$ is worked out in detail in \cite[Example 8.5.4]{egno}: again we recover that objects in ${\mathcal Z}(_{kG}{\mathcal M})$ are $G$-equivarient $G$-graded vector spaces.

\subsection{A geometric viewpoint on $_{D(G)}{\mathcal M}$}

Viewing the finite group $G$ as an algebraic or differential variety with finitely many points, with the natural $G$-action by conjugation, 
the category of  $D(G)$-modules is then seen in geometric terms, allowing for generalisations to the context of algebraic groups, Lie groups, and their analogues. We sketch here the different points of view.  

\subsubsection{Vector bundles on $G$}\label{sec:vect-bun}
Through geometric lens, a $G$-graded vector space  $V=\oplus_{g\in G} V_g$ with support ${\rm supp}(V)=\{h\in G~|~V_h\neq 0\}$ can be viewed as a vector bundle $\pi\colon \mathcal V=\coprod_{g\in G}V_g\to {\rm supp}(V)$, with $\pi(V_g)=g$. A $G$-equivariant vector bundle on $G$ is 
a vector bundle $({\mathcal V},\pi)$ on $G$ such that $\pi$ is $G$-equivariant, that is, $\pi(g.v)=g\pi(v_h)g^{-1}$ for all $g,h\in G$ and $v\in V_h$, \cite[Section 2.2]{arcata}.

This interpretation establishes an equivalence of categories between $_{D(G)}\mathcal M$ and the category of $G$-equivariant vector bundles on $G$, \cite[Theorem 2.2]{wither}. This interpretation lies behind Lusztig's theory of non-abelian Fourier transform that we will discuss in Section \ref{sec:Lusztig}. 

\medskip

Even if the terminology of ``center of a monoidal category'' is not explicitely mentioned, \cite[\S 8.4]{SS} contains a very detailed account on how to view the center of the monoidal category of $C(G)$-modules as the braided monoidal category of $G$-equivariant vector bundles: here the elements in $kG=C(G)^*$ are seen as measures on $G$. The authors declare to have learnt this from D. Kazhdan. 

\medskip

It reads as follows: left (or right) $C(G)$-modules are $kG$-comodules, i.e., vector bundles on $G$. The tensor product of two $G$-equivariant vector bundles $\mathcal V=\coprod_{g\in G}V_g$ and $\mathcal W=\coprod_{h\in G}W_g$ on $G$ has fiber $\bigoplus_{hl=g}V_h\otimes W_l$ on $g\in G$ and an element of the center is a pair $(\mathcal V,c_{-,{\mathcal V}})$ where $c_{-,{\mathcal V}}\colon {\mathcal W}\otimes {\mathcal V}\to {\mathcal V}\otimes {\mathcal W}$ preserves the $G$-grading for any vector bundle ${\mathcal W}$ on $G$. In particular, this must hold for the bundle ${\mathcal W}^h$ whose support is concentrated in $h\in G$  and with $1$-dimensional fiber $W_h$. Thus $c(h):=c_{\mathcal W^h,\mathcal V}$ on the fiber over $l\in G$ gives an isomorphism $V_{h^{-1}l}\to V_{lh^{-1}}$. Taking $l=hg$ we see that $c(h)$ gives an isomorphism between the fibers $V_g$ and $V_{hgh^{-1}}$. Now, $\mathcal W^h\otimes\mathcal W^l=\mathcal W^{hl}$ and so \eqref{eq:center-cat} shows that $c(hl)=c(h)c(l)$, for $h,\,l\in G$, that is, $(\mathcal V,c_{-,{\mathcal V}})$ is a $G$-equivariant vector bundle on $G$ or, equivalently, a $D(G)$-representation.
The condition required to morphisms in ${\mathcal Z}({\mathcal M}^{kG})$ is equivalent to the condition of $D(G)$-module morphisms.


\medskip

Observe that due to our choice of the presentation of $D(G)$, making use of $(kG^*)^{cop}$, the equivalence of categories induced from this identification does not identically map the tensor product of $D(G)$-modules to the tensor product of the corresponding vector bundles, but rather to the tensor product in the reversed order. However, the braiding in  $_{D(G)}\mathcal M$ can be used to obtain an equivalence of  monoidal categories, \cite[Exercise 8.1.9]{egno}.

\subsubsection{Sheaves on $G$}

Through algebro-geometric lens,  a $G$-grading on a vector space $V=\bigoplus_{g\in G}V_g$ determines a sheaf $\mathcal V$ on $G$: any subset $U$ of $G$ is open, the sections $\mathcal V(U)$ of $\mathcal V$ over $U$ are given by $\bigoplus_{g\in U}V_g$, and the restrictions  $\mathcal V(U)\to \mathcal V(U')$ for $U'\subset U$ are the projections along the decomposition given by the grading. Clearly this assignment satisfies the locality and gluing conditions. In particular, the stalk at $g\in G$ of $\mathcal V$ is the homogeneous component $V_g$. Just as for representations of algebras, the  tensor product of sheaves on a variety is not well-behaved and some caution is needed. A well-behaved operation is the external product of sheaves. For sheaves $\mathcal V$ and $\mathcal W$ on $G$, this is the sheaf on $G\times G$ given as follows: let $p_1,\,p_2\colon G\times G\to G$ be the canonical projections. These maps are open, so the pull-backs $p_1^*\mathcal V$ and $p_2^*\mathcal W$ are the sheaves on $G\times G$ whose stalks at $(g,h)$ are
$(p_1^*\mathcal V)_{(g,h)}=V_g$ and $(p_2^*\mathcal V)_{(g,h)}=V_h$, \cite[2.6.2]{springer}. The external tensor product $\mathcal V\boxtimes\mathcal W$ of $\mathcal V$ and $\mathcal W$ is then the sheaf on $G\times G$ whose stalk at the point $(g,h)$ is  $(p_1^*\mathcal V)_{(g,h)}\otimes (p_2^*\mathcal V)_{(g,h)}=V_g\otimes V_h$. In the special case in which the variety is a group, we may invoke the multiplication map $\mu\colon G\times G\to G$ and define also the (internal) tensor product $\mathcal V\otimes\mathcal W$ of $\mathcal V$ and $\mathcal W$ as the push forward $\mu_*(\mathcal V\boxtimes\mathcal W)$ along $\mu$ of $\mathcal V\boxtimes\mathcal W$. Its stalk at the point $g$ is then
$(\mathcal V\otimes\mathcal W)_g:=(\mu_*(\mathcal V\boxtimes\mathcal W))_g=(\mathcal V\boxtimes\mathcal W)(\mu^{-1}(g))=\oplus_{(h,l)\in G\times G, hl=g}V_h\otimes W_l$. In other words, the category of sheaves on $G$ and the category of $kG$-comodules are monoidally equivalent.   

\medskip

Equivariant sheaves are slightly harder to define than in the case of vector bundles: denote the action morphism given by conjugation of $G$ on itself by $\gamma\colon G\times G\to G$. Then a $G$-equivariant sheaf is the datum of a pair $(\mathcal V,\phi)$ where $\mathcal V$ is a sheaf on $G$ as above and $\phi$ is an isomorphism
of sheaves between the pull-backs $\gamma^*\mathcal V\to p_2^*\mathcal V$ satisfying the compatibility condition on $G\times G\times G$
\begin{equation}
(\mu\times\id_G)^*(\phi)=p_{23}^*(\phi)\circ(\id_G\times\gamma)^*(\phi)
\end{equation}
where $p_{23}\colon G\times G\times G\to G\times G$ denotes the projection to the second and third component.
Now, the stalk at the point $(g,h)$ of $\gamma^*(\mathcal V)$ is $V_{\gamma(g,h)}=V_{ghg^{-1}}$ hence at the level of stalks $\phi$ gives isomorphisms 
$\phi_{(g,h)}\colon V_{ghg^{-1}}\to V_h$. The compatibility condition at the stalk $(g,h,l)$ becomes 
$\phi_{(gh,l)}=\phi_{(h,l)}\circ \phi_{(g,hlh^{-1})}$. Setting $\rho(g):=\bigoplus_{h\in G}\phi_{(g,h)}^{-1}\colon \bigoplus_{h\in G}V_h\to  \bigoplus_{h\in G}V_{ghg^{-1}}$ we obtain a representation of $D(G)$ on $V=\bigoplus_{h\in G}V_h$. Morphisms in the category of $G$-equivariant sheaves correspond to morphisms of $D(G)$-modules. 

\medskip

Just as in the case of vector bundles, by our choice of realization of $D(G)$ the above identification induces an equivalence of categories reversing the monoidal structure. Combining with the braiding one obtains an equivalence of monoidal categories, \cite[Exercise 8.1.9]{egno}.

\subsection{Irreducible $D(G)$-modules}\label{classifica}

Assume now that ${\rm char} (k)$ does not divide $|G|$. Then, finite-dimensional $D(G)$-modules are semisimple. If ${\rm char}(k)=0$ this can be proved as in \cite{DPP}: in this case $D(G)$ is a based ring as defined in \cite[Section 1]{arcata}, and it is therefore semisimple.  More generally, semisimplicity can also be proved  using the existence of a left integral by standard Hopf algebra theory, as in \cite[Theorem 3.2]{gould}, where the case $k=\Cc$ was considered. 
All irreducible representations are finite-dimensional because $D(G)$ is finite-dimensional. 

The classification of irreducible $D(G)$-modules originally appeared in at least four different contexts independently:  in \cite[2.2]{DPP} and later in \cite[Section 6]{gould} in terms of the Drinfeld double itself, in \cite[Section 2.2]{arcata} in the language of $G$-equivariant vector bundles, and in \cite[Proposition 3.3]{CR} in the language of Hopf bimodules, where the authors compare the result to that of \cite{DPP}. We recall here the classification. In the successive sections we will review  the different historical motivations behind it. 

\bigskip
The decomposition \eqref{eq:deco-category} shows that each irreducible representation of $D(G)$ lies in a  subcategory $_{D(H)}{\mathcal M}^{\Oc}$ for some uniquely determined conjugacy class $\Oc$ in $G$. 
In addition, for each irreducible representation $V=\oplus_{h\in \Oc}V_h$, the homogeneous component $V_g$ is a representation of the centraliser $C_G(g)$, which is necessarily irreducible. Moreover, ${\rm Ind}_{C_G(g)}V_g=k G\otimes_{k C_G(g)}V_g\simeq V$, and $V_{xgx^{-1}}=k x\otimes_{k C_G(g)}V_g$. 
Conversely, all irreducible $D(G)$-representations are constructed this way. In fact, the correspondence between $C_G(g)$-representations and $D(G)$-representations by means of induction can be seen as a Morita equivalence, \cite[Section 2, Remark (III)]{mason}. 

Choosing a representative $g_{\Oc}\in G$ for each conjugacy class $\Oc$, the set of isomorphism classes of simple objects in $_{D(G)}{\mathcal M}^{\Oc}$ is then in bijection with the set of isomorphism classes of irreducible representations of $C_G(g_{\Oc})$ and the set of isomorphism classes of irreducible representations of $D(G)$ is in bijection with the set ${\mathcal M}(G)$  of $G$-orbits of pairs $(\Oc,\rho)$ where $\Oc$ is a conjugacy class in $G$ and $\rho$ is an irreducible representation of  $C_G(g_{\Oc})$.

\medskip

If ${\rm char}(k)$ divides $|G|$, then $D(G)$ is not semisimple but the set of isomorphism classes of irreducible modules is again in bijection with ${\mathcal M}(G)$ and the set of isomorphism classes of indecomposable modules is in bijection with the set of $G$-orbits of pairs $(\Oc,\rho)$ where $\Oc$ is a conjugacy class in $G$ and $\rho$ is an indecomposable representation of  $C_G(g_{\Oc})$, \cite[Proposition 1.2, Corollary 2.3]{wither}. 

\medskip

If $k=\mathbb C$ is observed in \cite{tom} that these representations preserve the $\star$-structure and form an exhaustive list of irreducible $\star$-representations of $D(G)$.

\medskip

It is of interest for category theorists to estimate the number of simple objects in a braided, or monoidal semisimple category. The following estimate of the number $f(n)$ of irreducible complex $D(G)$-modules depending on the number $n$ of irreducible representations of $G$ is due to P. Etingof, \cite[Theorem 3.1]{etingof} and answers a question of E. Rowell 
\begin{equation*}C_1n^{1/3} \log(n)\leq log f (n)\leq C_2n \log^7(n)\end{equation*}
for some $C_1,\,C_2\in{\mathbb R}_{>0}$.

\section{Occurrences of the Drinfeld Double}

In this section we single out a few further themes in mathematics and physics where the $D(G)$ or its representations play a role.

\subsection{Character theory, link invariants and generalized Thompson series}

Independently of \cite{arcata,DPP}, Gould has classified complex irreducible representations of $D(G)$ with the  purpose of producing  link invariants and  unitary representations of the braid group. Indeed, the action of the universal $R$-matrix on the two-fold tensor product of a module of a quasitriangular Hopf algebra $V$ produces a solution of \eqref{eq:braiding}. A recipe on how to produce link invariants from the $R$ matrix and $V$ is to be found in \cite[15.2]{CP}. 

With this motivation,  Gould computed  the matrix coefficients of irreducible $D(G)$-representations and the action of the $R$-matrix on tensor products, and developed the character theory. He remarked that the left integral used to prove semisimplicity of $D(G)$ is of use for constructing representations of the Temperley-Lieb algebra, and for retrieving $R$-matrices with spectral parameters \cite[Remark p. 285]{gould}.
He also observed that every irreducible representation of $D(G)$ is equivalent to a unitary one \cite[Lemma 4.1]{gould} generalising the averaging process for $\Cc G$. As a matter of fact, any finite-dimensional representation of $D(G)$ is equivalent to a unitary one. This is a consequence of the fact that  $D(G)$ is a finite-dimensional Hopf $C^*$-algebra, so its convolution $\star$-algebra is a finite-dimensional $C^*$-algebra. 

\medskip

One should mention that a complex character theory for a family of algebras including $D(G)$ has also been developed in \cite{Ba}, where it is suggested to use the characters of $D(G)$ to produce generalized Thompson series associated with $G$. The Thompson series for the Monster group were at the base of the Moonshine conjecture, establishing a connection between the characters of the Monster and modular forms. Formulated in \cite{CN} it was proved by Borcherds in \cite{B}. 

\medskip

For potential applications in conformal field theory, statistical mechanics, quantum integrable systems, and quantum computing, explicit matrix elements and character tables for $D(G)$ for $G$ ranging in several series of finite groups including cyclic and dihedral groups, and finite subgroups of ${\mathrm SU}(2)$ and ${\mathrm SU}(3)$ have been  calculated, \cite{dancer,Coq}.

\medskip

Back to the representations of braid groups, it is relevant to describe the image of  braid groups through representations stemming from the action of a quasitriangular Hopf algebra  $\mathcal H$. In particular,  it is of interest to determine under which conditions the image is finite. This has consequences also for the applications because the computational power of a topological quantum computation associated with these representations is strictly related to the size of the image, \cite[Examples 1.1, 1.2]{rowell2}. If $\mathcal H$ is a quantized enveloping algebra at a root of unity, then the image is finite only in special cases. However, it is shown in \cite{erw} that if $\mathcal H$ is the double of a finite group $G$, then the image is always finite. In addition, answering a question by Drinfeld, the authors show that if $G$ is a $p$-group, the image of the pure braid group is again a $p$-group. 

\subsection{Character theory of finite groups of Lie type}\label{sec:Lusztig}

The sets ${\mathcal M}(G)$ defined in Section \ref{classifica}  play an important role in Lusztig's classification of irreducible representations of  finite groups of Lie type. Roughly speaking, these groups are obtained from a reductive algebraic group $\mathbb L$ over the algebraic closure of a finite field $\mathbb F_q$, by taking points over the finite field itself, e.g., $L=\GL_n({\mathbb F}_q)\leq \GL_n(\overline{\mathbb F}_q)={\mathbb L}$. The key idea in the theory is that the geometry and the combinatorics related to the algebraic group $\mathbb L$ allow to give a description of the characters of $L$  that is uniform for all choices of $q$. This is obtained through a series of steps. First of all,  an inductive  procedure reduces the problem of the classification of characters to the case of a smaller class, that of unipotent ones, of possibly smaller groups that are again of Lie type. Unipotent characters in these groups are defined through the geometry of the corresponding algebraic group.     

\medskip

Next step is the study of the set $\mathcal U$ of irreducible unipotent characters  of a (split) finite group of Lie type $L$. The set $\mathcal U$ is parted into a disjoint union of subsets $\mathcal U_G$, where $G$ runs through a finite collection $\mathcal G$ of finite groups. These groups are the so-called canonical quotients:   they arise as quotients of component groups $C_{\mathbb L}(\mathrm u)/C_{\mathbb L}(\mathrm u)^\circ$ of suitably chosen unipotent elements $\mathrm u$ in $\mathbb L$. 

\medskip

The following step is to determine $\mathcal U_G$ for every $G\in\mathcal G$ 
and here is where the sets ${\mathcal M}(G)$ play a role, as ${\mathcal U}_G$ is in bijection with ${\mathcal M}(G)$, \cite[Section 4]{Orange}. The geometric interpretation of $D(G)$-modules as in Subsection \ref{sec:vect-bun} allows to identify the $\mathbb C$-linear span $\mathbb C\,\mathcal U_G$ of ${\mathcal U}_G$ with the complexified Grothendieck ring of $G$-equivariant vector bundles on $G$.

\medskip

This geometric insight then allows the determination of the characters in each subset ${\mathcal U}_G$ via the calculation of the so-called almost characters related to ${\mathcal M}(G)$. These are certain class functions on $L$ whose values are in principle computable as they can be interpreted as characteristic functions (up to a scalar) of character sheaves, \cite{Shoji}. Then the irreducible characters in ${\mathcal U}_G$ can be computed if one knows the matrix of the change of basis of  $\mathbb C\,\mathcal U_G$ between the basis of almost characters to the basis of unipotent irreducible characters $\mathcal U_G$. In \cite[4.14]{Orange} Lusztig defined a pairing that allows to derive almost characters as linear combinations of  irreducible unipotent characters, \cite[4.23]{Orange}, \cite[Theorem 14.2.3]{DM}. Even if the target of Lusztig's pairing is the algebraic closure $\overline{\mathbb{Q}}_\ell$ of the field of $\ell$-adic numbers, for $\ell$ a prime coprime with $q$, the pairing  takes in fact values in the subfield of the algebraic integers. Hence, it can also be regarded as a pairing with image in ${\mathbb C}$. It is defined on $\mathcal M(G)$ as follows:
\begin{align} \label{eq:Fourier}
\{ \  , \  \}&\colon \mathcal{M}(G)\times \mathcal{M}(G)\rightarrow \overline{\mathbb{Q}}_\ell\\
\{ ({\mathcal{O}_1}, \rho_1), ({\mathcal{O}_2}, \rho_2)\}&=\frac{|\mathcal O_1||\mathcal O_2|}{|G|^2}\sum_{\substack{h\in G\\ 
hg_{\mathcal{O}_2}h^{-1}\in C_G(g_{\mathcal{O}_1})}} Tr(\rho_1(h g_{\mathcal{O}_2}h^{-1})) Tr(\rho_2(h^{-1} g_{\mathcal{O}_1}h)),\nonumber
\end{align}
where as in Section \ref{classifica} the symbol $g_{\mathcal{O}}$ stands for a representative of the conjugacy class $\mathcal{O}$ in $G$. Identifying then ${\mathcal U}_G$ with ${\mathcal M}(G)$ and 
viewing a $\mathbb C$-linear combination of elements in $\mathcal U_G$ as a function on $\mathcal{M}(G)$, one obtains from \eqref{eq:Fourier}  the map
\begin{align}\label{eq:Fourier2}
FT_G\colon \mathbb C\mathcal U_G&\longrightarrow \mathbb C\mathcal U_G\\
f{(\Oc,\rho)}&\mapsto \sum_{({\Oc'},\rho')\in{\mathcal M}(G)}\{ ({\mathcal{O}}, \rho), ({\mathcal{O}'}, \rho')\}f{({\Oc'},\rho')}\nonumber
\end{align}
in analogy with the classical Fourier transform. 
The matrix associated with \eqref{eq:Fourier} has rows and columns indexed by $\mathcal{M}(G)$ and is called the non-abelian Fourier transform matrix associated with $G$. It is a unitary and Hermitian matrix. It may be seen as a normalized character table of $K_0(_{D(G)}\mathcal M)$, the Grothendieck group of ${D(G)}$, \cite[Section 14.2]{DM}.  Heuristically speaking, combining the non-abelian Fourier transforms for every $G\in \mathcal G$ gives a change of basis between $\mathcal U$ and the set of almost characters, with associated block diagonal matrix. In addition, the groups $G\in\mathcal G$ do not depend on $q$ \cite[Section 14.2]{DM} and are rather small, the biggest one being $\mathbb S_5$. Hence, the size of each block in the change of basis matrix is ``small''. The theory can be adapted to the case in which $L$ is not split, by suitably modifying the subsets $\mathcal M(G)$ and the transform. 

\medskip

The non-abelian Fourier transform for finite groups of Lie type has inspired the work on the representation theory of other families of groups. Notably, it has been used to define non-abelian Fourier transforms for unipotent representations of $p$-adic groups in a way that is compatible with reductions to those finite groups of Lie type that can be naturally obtained as subquotients od $p$-adic groups,  \cite{CO, ACR}.

\subsection{Hopf algebra structures on path algebras and quantum groups}

In the quest for examples of non-commutative non-cocommutative Hopf algebras, Cibils and Rosso considered possible Hopf algebra structures on the class of  path algebras of a quiver. The class of path algebras is not a niche choice: if $k$ is algebraically closed then the category of representations of any finite-dimensional sufficiently well-behaved algebra (i.e., hereditary), is equivalent to the category of representations of a finite-dimensional path algebra, \cite[I, Corollary 6.10, VII. Theorem 1.7]{ass}.  

We recall that a {\em quiver} $\mathcal Q=(\mathrm V,\mathcal A)$ is a directed graph where loops are allowed: $\mathrm V$ denotes the set of vertices and  $\mathcal A$ the set of arrows. A path of length $r>0$ in $\mathcal Q$ is then a sequence $a_r\cdots a_1$ with $a_i\in \mathcal A$ for every $i\in\{1,\,\ldots,r\}$, such that the target of $a_i$ equals the source of $a_{i+1}$ for every $i\leq r-1$. Each vertex is considered to be a path of length $0$. The source (respectively target) of a path is the source  (respectively target) of its rightmost (respectively leftmost) term. The {\em path algebra} $k\mathcal Q$  of $\mathcal Q$ has the formal $k$-span of paths as underlying vector space, and product given by concatenation, with the understanding that for two paths $p_1$ and $p_2$, there holds $p_1p_2=0$ if the source of $p_1$ differs from the target of $p_2$. It is naturally graded by the length of the paths, non-commutative for any non-trivial $\mathcal Q$, and finite dimensional if and only if $\mathrm V$ and $\mathcal A$ are finite sets and $\mathcal Q$ has no oriented cycles. 

Cibils and Rosso observed that if $k\mathcal Q$ admits a Hopf algebra structure, then: the algebra $k^{\mathrm V}$ of $k$-valued functions on $\mathrm V$ with pointwise multiplication admits a Hopf algebra structure;  
the $k$-span $k\mathcal A$ of $\mathcal A$ is a Hopf bimodule for $k^{\mathrm V}$; and $k\mathcal Q$ is obtained as a $k^{\mathrm V}$-tensor algebra of $k\mathcal A$ over $k^{\mathrm V}$. The Hopf algebra structures on $k^{\mathrm V}$ are in bijection with the group structures on $\mathrm V$, that is, $k^{\mathrm V}\simeq (k\mathrm V)^*$ for some group  structure on $\mathrm V$, and the problem boils down to understanding the Hopf bimodule structures on $(k\mathrm V)^*$ for any group structure on $\mathrm V$. In virtue of the equivalences in Section \ref{sec:bimodules} and the isomorphism \eqref{eq:dual}, this corresponds to understanding  $D(\mathrm V)$-modules and then translating this knowledge into recipes for assigning the arrows. This is done by means of a variant of the Cayley graph for the group $\mathrm V$. Specifically, the dimension of the $\mathrm V$-graded components of a $D(\mathrm V)$-module gives a function $m\colon \mathrm V\to\mathbb N_{\geq1}$ that is constant on conjugacy classes of $\mathrm V$: the quiver $Q=(\mathrm V,\mathcal A)$ obtained assigning $m(g_1^{-1}g_2)$ arrows from $g_1$ to $g_2\in \mathrm V$  has a path algebra admitting a graded Hopf algebra structure, and all graded Hopf algebra structures on a path algebra are obtained in this way,  \cite[Th\'eor\`eme 3.1]{CR}. The authors do not restrict themselves to the case of finite groups, and apply their result to realize the  quantized enveloping algebra of a semisimple Lie algebra with symmetric Cartan matrix as a Hopf algebra quotient of the path algebra of a suitable quiver whose set of vertices is an abelian group. 

\subsection{Classification of Hopf algebras, Nichols algebras and quantum groups}
The representations of $D(G)$ are crucial also in Andruskiewitsch and Schneider's lifting method for classifying finite-dimensional (pointed) Hopf algebras, \cite{AS}. We give here a brief account of their method, assuming that $k$ is algebraically closed and of characteristic zero. 

\medskip

To any Hopf algebra ${\mathcal H}$ one may attach two natural invariants: its {\em co-radical} ${\mathcal H}_0$, i.e., the sum of its simple subcoalgebras, and the group $G(\mathcal H)=\{h\in \mathcal H~|~\Delta(h)=h\otimes h\}$ of its group-like elements. Clearly, $kG(\mathcal H)\subseteq {\mathcal H}_0$,  
if equality holds, then ${\mathcal H}$ is called {\em pointed}. The class of pointed Hopf algebras includes several important examples: group algebras (here $\mathcal H=kG(\mathcal H)=\mathcal H_0$), universal enveloping algebras of semisimple Lie algebras (here $G(\mathcal H)=1$), and their quantum analogues (here $G(\mathcal H)$ is the subalgebra generated by the Cartan part). Now, the coradical induces a natural filtration of $\mathcal H$ and the associated graded Hopf algebra is a Radford biproduct (a construction similar to the semi-direct product) of $kG(\mathcal H)$ and an $\mathbb N$-graded subalgebra $\mathcal N$, on which $kG(\mathcal H)$ acts by graded automorphisms. The degree $1$ component $\mathcal N_1$ of $\mathcal N$ has a natural action of $D(G(\mathcal H))$ and the relations in the subalgebra generated by $\mathcal N_1$ are completely determined by the braiding stemming from the $D(G(\mathcal H))$-action. It is the graded (braided) bialgebra that one naturally associates to a vector space with a braiding by means of the procedure introduced in \cite{nichols} and re-discovered in \cite[Section 3]{wor-diff-cal} as algebras of coinvariants of quantum differential forms. Such a bialgebra is nowadays called a {\em Nichols algebra}, or {\em Nichols-Woronowicz algebra}, or {\em small shuffle algebra};  it was called a {\em bitensor algebra} in \cite{Sch2}.

\medskip

Concretely, the construction of the Nichols algebra for a vector space $V$ with braiding $c\in {\rm Aut}(V^{\otimes 2})$ reads as follows. For any $n\geq 2$ and any $\sigma\in\mathbb S_n$,  choose a reduced expression $s_{i_1}\cdots s_{i_r}$ of $\sigma$ as a product of transpositions $s_i=(i,\,i+1)$ for $i\in\{1,\,\ldots,\,n-1\}$, and consider the operator $M_n(\sigma)\in{\rm Aut}(V^{\otimes n})$ given by 
$c_{i_1}\cdots c_{i_r}$ where $c_{i}\in{\rm Aut}(V^{\otimes n})$ acts as $c$ on the $i$ and $(i+1)$-th tensor factor and trivially on the others.  By \cite{Matsu} the operator $M(\sigma)$ does not depend on the choice of the reduced decomposition. As an algebra, the Nichols algebra is the quotient of $T(V)$ by the kernel of the  {\em quantum symmetrizer} $\bigoplus_{n\geq 2}\sum_{\sigma\in\mathbb S_n}M_n(\sigma)$. The coalgebra structure is determined by taking all elements in $V$ to be primitive. 

\medskip

The family of Nichols algebras includes many different and notable examples: the symmetric algebra of a vector space $V$ is the Nichols algebra corresponding to the braiding $\tau$; the exterior algebra of $V$ is the Nichols algebra corresponding to the braiding $-\tau$; the positive part of the quantized enveloping algebra of a semisimple Lie algebra $\mathfrak g$ is the Nichols algebra corresponding to the span of simple root vectors $E_1,
,\ldots E_n$ with braiding $c(E_i\otimes E_j)=q^{d_{ii}a_{ij}}E_j\otimes E_i$, where $A=(a_{ij})$ is the Cartan matrix of $\mathfrak g$ and $D=(d_{ij})$ is a symmetrising matrix for $A$ \cite{rosso,rosso-inv, Sch2,Mu}. In other words,  the algebra $\mathbf f$ in \cite[Chapter 3]{Lu-QG} is a Nichols algebra.  

\medskip

Back to the pointed Hopf algebras classification problem, one can fix a finite group $G$ and determine all pointed Hopf algebras $\mathcal H$ such that $G(\mathcal H)=G$ according to the following strategy: 
single out all representations of $D(G)$ for which the corresponding Nichols algebra is finite-dimensional; write a presentation for each of these Nichols algebras; determine all finite-dimensional  Hopf algebras $\mathcal H$ such that their associated graded Hopf algebra is the biproduct of $kG$ and any of the algebras in the list; verify if all finite-dimensional pointed Hopf algebras with group of grouplikes isomorphic to $G$ are indeed generated in degree $0$ and $1$. 

\medskip

For abelian groups this program has been completed by Andruskiewitsch, Angiono, Heckenberger, Schneider and their collaborators, using Lie theoretical methods. It showed striking connections with the classification of semisimple Lie algebras and Lie superalgebras in positive characteristic, and with the description of quantized enveloping algebras of Lie algebras and Lie superalgebras for which the parameter is a root of unity. For an account of the theory see \cite{AA,HS} and references therein.  

\medskip

When $G$ is non-abelian the classification is more involved and the program has been carried out only for very few examples. Finite-dimensional Nichols algebras associated with braidings coming from $D(G)$ for $G$ non-abelian seem to be very scarce. 
A big step forward in the case in which $G$ is solvable has been recently taken in \cite{AHV}, which contains, among other results, the classification of finite-dimensional pointed Hopf algebras of odd dimension. 

\medskip

A folklore conjecture states that if $G$ is simple and non-abelian there is no finite-dimensional Nichols algebra associated with a $D(G)$-module, implying that $\mathbb C G$ is the only finite dimensional complex pointed Hopf algebra with group of grouplikes isomorphic to $G$. The conjecture has been confirmed  for the alternating groups, \cite{AFGV1}, for all sporadic groups except from the Fischer group $Fi_{22}$,  the Baby Monster and the Monster, \cite{AFGV2}, and for several families of finite simple groups of Lie type, see \cite{CC, ACG5, ACG7}, but the general picture is still being figured out.

\subsection{Mapping class groups of surfaces and ${\mathrm SL}_2(\mathbb Z)$-actions}\label{sec:mapping}

We recall from \cite{RT} that a quasitriangular Hopf algebra $\mathcal H$ with universal $R$-matrix $\mathcal R$, monodromy element $\mathcal Q$ and Drinfeld element $\mathrm u$ is called a {\em ribbon Hopf algebra} if it has a central element $\mathrm v$, called {\em ribbon element}, satisfying the relations 
\begin{equation}
\Delta_{\mathcal H}(\mathrm v)=\mathcal Q^{-1}(\mathrm v\otimes\mathrm  v),\quad S_{\mathcal H}(\mathrm v)=\mathrm v,\quad \mathrm  v^2=\mathrm  uS(\mathrm u), \quad \varepsilon_{\mathcal H}(\mathrm v)=1.
\end{equation}
If this is the case, the action of $\mathrm v$ on any $\mathcal H$-module equips $_{\mathcal H}\mathcal M$ of a ribbon category structure, \cite[Chapter XIV]{kassel},  a notion introduced by Joyal and Street \cite{JS}, that is useful to provide isotopy invariants of framed links.

\medskip

The Drinfeld double of a Hopf algebra $H$ need not be a ribbon Hopf algebra, but if  $H=D(G)$ this is always the case.  Indeed,  for the $R$-matrix $R_{D(G)}$ (respectively,  $R_{D(G)}'$) as in \eqref{eq:errematrice} the monodromy element is 
$Q=\sum_{h,g\in G}\delta_hg\otimes\delta_{g^{-1}}g^{-1}h^{-1}g$ 
(respectively, $Q'=\sum_{h,g\in G}\delta_{hg^{-1}h^{-1}}h\otimes\delta_{h^{-1}}g$). The Drinfeld element  $u=\sum_{g\in G}\delta_{g^{-1}}g=S_{D(G)}(u)$ (respectively $u'=\sum_{g\in G}\delta_gg$) is central, so $v=u$  (respectively, $v'=u'$) is a ribbon element. Through the identification \eqref{eq:identify}, the ribbon element $u$ (respectively, $u'$) becomes the linear extension of the characteristic function of the diagonal (respectively, of the graph of the inversion map) in $G\times G$. 

\medskip

The existence of a ribbon element is very useful in topology and specifically in the study of the {\em mapping class groups} of surfaces. These are the groups of isotopy classes of homeomorphisms of surfaces, possibly preserving some marked points and/or orientation. Any factorizable ribbon Hopf algebra $\mathcal H$ with monodromy element $\mathcal Q=\mathcal Q^{(1)}\otimes \mathcal Q^{(2)}$ and ribbon element $\mathrm v$ gives rise to a (projective) representation of the mapping class group of any compact oriented surface of genus $g$ and $n$ marked boundary circles \cite{ly}, see also \cite[Section 5]{tu}. The representation space is ${\rm Hom}_{\mathcal H}(X^{\otimes n},(\mathcal H^*)^{\otimes g})$ where $X$ is an object in $_{\mathcal H}\mathcal M$ and $\mathcal H^*$ is equipped with the coadjoint action \eqref{eq:coadjoint}. Just as for the image of the braid group representations arising from the $R$-matrix, it is relevant to estimate the size of the image of these representations. An important feature of the representations arising from $D(G)$, for surfaces with at most one boundary point, is that the image is finite \cite{erw,FF}, see also \cite{gust} where the finite-image result is proved in the more general case of a twisted double and oriented compact manifolds with boundary, by different methods. 

\medskip

We have now a closer look at the case the torus, i.e., $g=1$ and $n=0$, following the exposition in \cite[Section 2]{FF}. Here, the mapping class group and the orientation preserving mapping class group act faithfully on the fundamental group of the torus $\mathbb Z\times\mathbb Z$ and are isomorphic, respectively, to $\mathrm{GL}_2(\mathbb Z)$ and $\mathrm{SL}_2(\mathbb Z)$. 
We recall that the latter is generated by the matrices $s=\left(\begin{smallmatrix}0&1\\-1&0\end{smallmatrix}\right)$ and $t=\left(\begin{smallmatrix}1&1\\0&1\end{smallmatrix}\right)$ and is completely determined by the relations $s^4=1$ and $s^2=(st)^3$. 

In this case, the representation space is 
\begin{equation*}
{\rm Hom}_{\mathcal H}(k,\mathcal H^*)\simeq (\mathcal H^*)^{\mathcal H}=\{f\in\mathcal H^*~|~{\rm Ad}^*(h)(f)=\varepsilon_{\mathcal H}(h)f,\;\forall h\in\mathcal H\}.
\end{equation*}
The $\mathrm{SL}_2(\mathbb Z)$-action is given by 
\begin{equation}\label{eq:st}
t\cdot f(-)=f(S_{\mathcal H}(\mathrm v^{-1})\,-), \mbox{ and }s\cdot f(-)=f(\mathcal Q^{(1)})\Lambda(\mathcal Q^{(2)}\,-)\end{equation} 
where $\Lambda\in \mathcal H^*$ is an integral, that is, it satisfies $f\Lambda=\varepsilon_{\mathcal H^*}(f)\Lambda$ for $f\in\mathcal H^*$.  Observe that since $\mathcal H$ acts on $\mathcal H^*$ by algebra automorphisms, $(\mathcal H^*)^{\mathcal H}$ is a subalgebra of $\mathcal H^*$.

Let us now specialize to the case in which ${\rm char} (k)$ is coprime with $|G|$ and  $\mathcal H=D(G)$. With the realization from Section \ref{sec:analyst} and making use of the non-degenerate pairing on $C(G\times G)$ given by
\begin{equation}\label{pairing}
    \langle f_1,f_2\rangle:=\int_G\int_Gf_1(g_1,g_2)f_2(g_1,g_2)dg_1dg_2,
\end{equation}
the vector space $D(G)^*$ is identified with $D(G)$, whence with $C(G\times G)$.
With this identification, multiplication, unit, comultiplicaton, counit and antipode of $D(G)^*$ become respectively
\begin{align*}
&(f_1\bullet f_2)(h,g)=|G|\int_G f_1( g_1,g)f_2(g_1^{-1}h,g)dg_1,&&1_{D(G)^*}(h,g)=\delta_{1_G}(h),\\
&(\Delta_{D(G)^*}(f))(h_1,g_1,h_2,g_2)=f(h_1,g_1g_2)\delta_{h_2, g_1^{-1}h_1 g_1},&&\varepsilon_{D(G)^*}(f)=|G|\int_Gf(g, 1_G)dg,\\
&S_{D(G)}^*(f)(h,g)=f(g^{-1}h^{-1}g,g^{-1}),
\end{align*}
for $h,h_1,h_2,g,g_1,g_2\in G$ and $f, f_1, f_2\in C(G\times G)$, see \cite{tom}. We spell out the coadjoint action of $D(G)$ on $D(G)^*$. For $g\in G$, the elements $\sum_{h\in G}\delta_{(h,g)}$ and $\delta_{(g,1_G)}$ generate $D(G)$. Their action on a function $f\in C(G\times G)\simeq D(G)^*$ is given by
\begin{align*}
(Ad^*(\sum_{h_1\in G}\delta_{(h_1,g_1)})f)(h_2,g_2)&=f(g_1^{-1}h_2g_1,g_1^{-1}g_2g_1),\\
(Ad^*(\delta_{(h_1,1_G)})f)(h_2,g_2)&=\delta_{[g_2^{-1},h_2^{-1}],h_1}f(h_2,g_2)
\end{align*}
for $g_1,g_2, h_1,h_2\in G$.
Therefore $(D(G)^*)^{D(G)}$ consists precisely of the functions that are supported on $(G\times G)_{comm}=\{(h,g)\in G\times G~|~hg=gh\}$ and are invariant under simultaneous $G$-conjugation of the variables. We denote the space of such functions by $C((G\times G)_{comm})^G$. The integral $\Lambda\in D(G)^*$ is uniquely determined up to a scalar and can be chosen to be the Haar functional $\widetilde{h}\colon D(G)\to {\mathbb C}$ given by $\widetilde{h}(f)=|G|\int_Gf(g,1_G)dg$. 
Then, the actions $\cdot$ and $\cdot'$, for the ribbon structures given by $R_{D(G)}, v$ and $R'_{D(G)}, v'$, respectively, are 
\begin{align}\label{eq:action}
(s\cdot f)(h,g)=f(g, h^{-1}),  &&(t\cdot f)(h,g)=f(h,h^{-1}g);\\
(s\cdot' f)(h,g)=f(g^{-1},h), &&(t\cdot' f)(h,g)=f(h,hg) &&
\mbox{ for }f\in C(G\times G)\mbox{ and }h,g\in G. 
\nonumber
\end{align}
They are obtained from one another  twisting by the automorphism of $\mathrm{SL}_2(\mathbb Z)$ that inverts $s$ and $t$. 
The $\cdot'$-action stems from the right action on $G\times G$ given by $(h,g)\left(\begin{smallmatrix}a&b\\
c&d\end{smallmatrix}\right)=(h^ag^c,h^bg^d)$ and for $A\in \mathrm{SL}_2(\mathbb Z)$ and $f\in C(G\times G)$ there holds $A\cdot' f(h,g)=((h,g)A)$. This action has been studied in \cite{tom}, where it was already observed that the $t$-action stems from the ribbon element and that the operators are well defined on the whole $C(G\times G)$, although they do not preserve the relation $(st)^3=s^2$, compare also with \cite[Section 4.b]{Dijk}. 

\medskip

We can say a bit more. The $\mathrm{SL}_2(\mathbb Z)$-actions on $C((G\times G)_{comm})^G$ naturally extend to a $\mathrm{GL}_2(\mathbb Z)$-action. We point out that the $\cdot'$-action of the involutions 
\begin{equation}\label{eq:j1j2}j_1:=\begin{pmatrix} -1&0\\0&1\end{pmatrix} \quad \mbox{ and }\quad j_2:=\begin{pmatrix} 1&0\\0&-1\end{pmatrix}\end{equation}
can be described in Hopf algebraic terms as 
\begin{equation}\label{eq:antipode-involution}
j_1\cdot' f:=f\circ (S_{kG^*}\otimes \id),\quad\mbox{ and }\quad  j_2\cdot' f:=f\circ (\id\otimes S_{kG}).   
\end{equation}
Even though $\id\otimes S_{kG}$ and $S_{kG^*}\otimes \id$ are not automorphisms of $D(G)$, one can directly verify that $j_i\cdot (f_1\bullet f_2)=(j_i(f_1)\bullet j_i(f_2))$ for $i=1,2$ whenever $f_i\in C((G\times G)_{comm})^G$. 

\bigskip

We point out that for  the exponent $e_G$ of $G$, an action of $\mathrm{GL}_2(\mathbb Z/e_G\mathbb Z)$ on $C((G\times G)_{comm})^G$ is defined in \cite[Sections 3.5.3, 8.4.2]{br}, \cite[18.5.1]{Br}, building on ideas from \cite[Ch. VII, \S 3]{dm}, in the framework of representations of finite groups of Lie type.

\subsection{Quantum Observables and Mackey quantization}

 When considering a physical system given as a particle living on a space $X$ and moving according to the action of a group $G$ on this space (and as such constrained on orbits) there is a  very general framework, known as \emph{Mackey quantization scheme} according to which the crossed product algebra $C(X) \rtimes G$ may be interpreted as a universal algebraic quantization of the system. 

In a topological setting this applies to the case of a topological group action of $G$ on $X$ (under suitable regularity assumptions) as being quantized by the (non commutative) transformation group $C^*$-algebra $C^*(G,C)=G\rtimes C_0(X)$ which plays the role of quantum observables; its unitary irreducible representations are the superselection sectors of the theory. 

When we restrict to homogeneous actions and $X$ can be identified with a group quotient $G/K$ with respect to the stabilizer $K$ of a point, these unitary irreducible representations can be constructed as induced representations starting from representations of the closed subgroup $K$. The Hilbert space on which the representation takes place is the space of $L^2$-sections of the vector bundle associated to the principal bundle $G\to G/K$. Mackey imprimitiviy theorem guarantees that all representations of interest can be constructed in this way.

The whole construction, however, may be rephrased in purely Hopf-algebraic terms and, as such, 
 this scheme is applied in \cite{maj-qm} to the case of the standard quantum algebra $U_q(su(2))$.

If $G$ is a finite group and $k=\mathbb C$, the algebraic part of the quantization scheme allows to understand $D(G)$ as the quantization of a particle constrained to move on a conjugacy class of $G$, as explained in \cite[Example 6.1.8]{mabook}.


\medskip

Recently Majid and McCormack \cite{majmc} proposed a different approach to irreducible $\star$-representations of $D(G)$ in the spirit of non commutative geometry.
Their aim is to mimick the geometric construction by Wigner of unitary irreducible representations of a semidirect product of Lie groups
$P=T\rtimes S$ (where tipically $T$ is assumed to be abelian and $S$ semisimple) on the space of sections of a suitable homogeneous vector bundle. 

Rather than looking for irreducible $\star$-representations of $D(G)$, in fact, they look after the equivalent set of irreducible 
$\star$-corepresentations of the dual algebra $D(G)^*$, thought of as the algebra of functions on the group-like object $D(G)$, where its semidirect product decomposition is reflected in the bicrossed product construction of $D(G)^*\simeq \mathbb CG\bowtie C(G)$, \cite[Example 7.2.5]{Maj}. This set of comodules is again in bijection with ${\mathcal M}(G)$.

Wigner construction requires, as initial data, a point $p_m$ in the  $P$-space $\widehat T$ together with an irreducible $\star$-representation 
of its stabilizer in $S$ (the little subgroup $L_m$ of $p_m$). The subindex $m$ parametrizes the set of orbits and can be understood
as the mass of a particle undergoing a quantization procedure. Together with the additional data of an irreducible $\star$-representation
of $L_m$ on a vector space $V_\pi$ this allows to construct an induced representation of $P$ on sections of the vector bundle
associated to the homogeneous $P$-space $P/(T\rtimes L_m)\simeq S/L_m$.  The image of the group elements as
operators on the (Hilbert) space of sections of the vector bundle are the quantum observables associated to the particle.

Here, replacing $P$ by the non commutative Hopf algebra $D(G)^*$, the decomposition $D(G)^*= \mathbb CG \bowtie  C(G)$ takes the place
of its semidirect product decomposition, with $\mathbb C G$ having the role of the abelian algebra. The starting point is, then, a $1$-dimensional corepresentation
of the group algebra, identified by an element $g_{\mathcal O}$ in a fixed conjugacy class $\mathcal O\subseteq G$. Its stabilizer $G_{\mathcal O}$ can be seen as a \emph{quantum subgroup}, i.e.
naturally defines a Hopf algebra surjection
$$ 
 C(G)\to  C(G_{\mathcal O}).
$$

In place of the $P$-space $P/(T\rtimes L_m)$ Majid and McCormack consider the Hopf-Galois 
extension corresponding to the projection

$$ 
D(G)^* \to  C(C_G)\bowtie \mathbb CG 
$$
 having a space of coinvariants isomorphic to $C(G/G_{\mathcal O})\simeq  C(\mathcal O)$. The Wigner induction procedure, then,
starts with an irreducible corepresentation of $ C(G_{\mathcal O})$ on a vector space $V_\pi$ which allows to construct
a corresponding cotensor product as a $ C(G/G_{\mathcal O})$-bimodule of coinvariants inside $D(G)^*\otimes V_\pi$.
Then, as shown in \cite[Lemma 3.3]{majmc}, the cotensor product $E$ is a $\star$-irreducible $D(G)^*$-comodule and every such
comodule can be constructed in this way, thus recovering the classification of $\star$-irreducible representations of $D(G)$ recalled in Section \ref{classifica}.

In the Wigner construction, averaging over the semidirect component $S$ with respect to a Haar measure allows to realize the induced
$\star$-representation on a subspace of sections of a vector bundle associated to the homogeneous principal bundle $P\to P/S$ in such a way that
the choice of the initial point $p_m$ is equivalent, after a suitable Fourier-type transform, to a Klein-Gordon equation for sections of the vector
bundle (after the usual identification of the quadratic Casimir of the universal enveloping algebra with invariant differential operators).

Majid and McCormack, therefore, have all the needed algebraic machinery at hand to perform analogous computations in the Drinfeld double case
starting from the Haar measures, in a couple of different directions. On one hand they can pass from the cotensor product $E$ to 
a left $D(G)^*$-subcomodule of $\mathbb CG\otimes W$, where $W$ is a right $D(G)^*$-comodule containing $V_\pi$ as subcomodule:
this in analogy with a Fourier transform where $ C(\mathcal O)$ has the role of momentum space and $\mathbb C G$ the role of position space. 
On the other hand, relying on the rich algebraic structure of the double they can rather start with the Hopf algebra surjection
$D(G)^* \to \mathbb C G$ to construct a dual interpretation of the transform as an explicit map embedding $D(G)^* \otimes V$, for any right $\mathbb CG$-comodule
$V$, inside $\mathbb C G\otimes W$ for a right $D(G)^*$-comodule $W$, under the condition that $V$ and $W$ have the same underlying $G$-grading.

To recover the interpretation of irreducible $\star$-corepresentations of $D(G)^*$ as invariant subspaces under suitable differential operators one needs to choose a suitable
bicovariant non commutative differential calculus on the quantum principal bundle we are working with.   However, using the fact that
$D(G)^*$ is coquasi-triangular and factorizable its bicovariant calculi are in $1-1$ correspondence with two-sided ideals of $D(G)$ as proven in \cite{maj-LNPAM},
thus ultimately again with irreducible representations of $D(G)$. Therefore any such calculus has the form
$$
D(G)^* \otimes (\bigoplus_{(\mathcal O,\pi)} \mathrm{End}(V_{(\mathcal O,\pi)}
))$$    
for a certain selection of pairs $(\mathcal O,\pi)$ of a conjugacy class and a corepresentation of its algebra of functions, inducing a $D(G)$-irreducible representation. When moving to quantum homogeneous spaces
one needs to consider $D(G)^*$ covariant calculi on $\mathbb C G$, and these correspond to two-sided ideals of $\mathbb C G^+={\rm Ker}(\varepsilon)$.
Once such a choice is made it is possible, much in the same way as done for quantum homogeneous spaces already in \cite{ditom}, to define quantum Laplacians over $\mathbb CG$ and study their eigenspaces
which allows to define the space carrying the irreducible corepresentation of $D(G)^*$ as a subspace of $\mathbb C G\otimes W$ satisfying a 
Klein-Gordon type equation (see Corollary 4.17 in \cite{majmc}).   

\subsection{Orbifold conformal field theory}


As we will see in Section \ref{sec:Fourier} the work of Tom Koornwinder and collaborators \cite{tom} was originally stimulated by the appearance of the quantum double of a finite group
in orbifold rational conformal field theory (RCFT) and the desire to give a purely representation theoretic explanation 
of the Verlinde formula for fusion rules. In fact, rational conformal field theory  provides  a somewhat  unifying bridge
giving sense to, if not explaining completely, the ubiquity of $D(G)$. In RCFT the Hilbert space of fields decomposes into a 
sum of representations of the left and right Virasoro algebra, called subsectors. Such space comes equipped with a 
vertex operator algebra  $\cal V$. In the holomorphic case $\cal V$ turns out to be also the only simple $\mathcal V$-module. 
The orbifold version is determined by the action of a  finite group  $G$ by automorphisms of the vertex operator algebra,
such that the invariant part ${\cal V}^G$ is completely reducible and has only a finite number of isomorphism classes of 
simple modules. Overall this guarantees the existence of a modular fusion category $\mathcal A_G$ that defines a 
fusion ring.

In \cite{Dijk} it was shown how this modular fusion category can be reconstructed, as a Tannaka-Krein type of duality, from the Drinfeld
double of the group $G$, and irreducible $D(G)$-modules were classified in order to 
recover the fusion rules, the S matrix and the conformal weights of an orbifold CFT. This approach was later generalized in \cite{DPP} to orbifold conformal field theories depending also on a group cocycle $\omega$,
showing that the fusion category is then isomorphic to the category of unitary irreducible representations of a twisted Drinfeld double.

On the other hand conformal field theories can also be described in more geometrical terms, as in \cite{Segal}, as associated to a Riemann surface with punctures
where to each puncture is associated the representation space of, say, a Virasoro algebra. Each such space can be decomposed into spaces of holomorphic sections of suitable flat vector bundles, which are $G$-equivariant in the orbifold case.  Fusion rules are then corresponding to 
the way in which data at punctures interact when sewing together such surfaces. Thus the geometrical and algebraic properties of orbifold RCFT provide, in principle, a general setting linked to group representations, mapping class group and the vector bundle interpretation. 

An account written for group theorists on how $D(G)$ occurs in conformal field theory and the theory of holomorphic orbifolds is \cite{mason}.

\subsection{Fourier transforms for $D(G)$ and Verlinde formulas}\label{sec:Fourier}

Beside Lusztig's non-abelian Fourier transform described in Section \ref{sec:Lusztig}, other variants of Fourier transforms for finite-dimensional Hopf algebras in general or specifically for $D(G)$ have been introduced, the first being probably in \cite{KP} in the framework of finite-dimensional $C^*$-algebras and Kac algebras. Notice that finite-dimensional Kac algebras are precisely the  finite quantum groups according to \cite{W}, so $D(G)$ fits into this context. In the setup of finite-dimensional semisimple Hopf algebras, a Fourier transform was defined in \cite{AN1,AN2} as a morphism of right $\mathcal H^*$-modules from a Hopf algebra $\mathcal H$ to its dual $\mathcal H^*$: for $\mathcal H$ being either $D(G)$ or $D(G)^*$ respectively, it corresponds in the notation from Section \ref{sec:analyst} and \ref{sec:mapping} to the transforms $F_{D(G)}(f)(h,g)=f(h^{-1},g)$ and 
$F_{D(G)^*}(f)(h,g)=f(g^{-1}hg,g^{-1})$ for $h,g\in G$ and $f\in D(G)\simeq C(G\times G)$. The authors also introduce a right-handed variant $\mathbb G_{\mathcal H}$, a notion of Hecke algebra for a Hopf algebra inclusion $K\to \mathcal H$, study spherical functions, and  define a notion of Gelfand pair in a Hopf algebra (see also \cite{Va}), showing that $(D(G),G)$ is always a Gelfand pair. A Fourier transform  for a factorizable ribbon Hopf algebra $\mathcal H$ with ribbon element $\mathrm v$ is given in \cite{LM}. This is an invertible linear map $S_-\colon \mathcal H\to\mathcal H$ satisfying $S_-^2=S_{\mathcal H}^{-1}$ and $(\lambda(\mathrm v))S_-^2= (S_-T)^3$, where $T$ is a linear map playing the role of an exponentiated Laplacian and $\lambda\in\mathcal H^*$ is a normalized integral. For $\mathcal H=D(G)\simeq C(G\times G)$ this is given by $S_-(f)(h,g)=f(hg^{-1}h^{-1},h)$ for $h,g\in G$ and $f\in C(G\times G)$. Also \cite{CW}  contains an analysis of a Fourier transform $\Phi$ for arbitrary finite-dimensional Hopf algebras. It is obtained by pre-composing the $\mathbb G_\mathcal H$-transform in \cite{AN1} with the antipode, and it is shown that for ribbon Hopf algebras the transform $S_-$ from \cite{LM} equals $f_Q\circ \Phi$, for $f_Q$ is as in \eqref{eq:fQ}. A definition of Fourier transform for arbitrary Hopf algebra is also present in \cite{kop} where it is used to derive an analogue of Heisenberg uncertainty principle for Hopf algebras through dimension estimates. 

\bigskip

 The $\cdot'$-action \eqref{eq:action} of the  element $s^{-1}\in{\mathrm SL}_2(\mathbb Z)$ on $C(G\times G)^G_{comm}$ has been interpreted as a Fourier transform for $D(G)$ in \cite{tom}, with the purpose of clarifying the Verlinde formula for the fusion rules for $D(G)$ from \cite{Dijk} in a purely algebraic setup without referring to conformal field theory. Relations with the action of $s$ emerge also for the other Fourier transforms, \cite{LM,CW,exotic}.  We sketch here the arguments in \cite{tom}, observing that most of them are effective whenever the characteristic of $k$ does not divide $|G|$.
 To lighten up notation, we omit the symbol $\cdot'$ to indicate the action. 

 \medskip
 
 As observed in \cite[p. 289]{gould}, the properties of the traces ensure that the characters of $D(G)$ live in the subalgebra  of $D(G)^*$ consisting of invariants for the coadjoint action, that is $C((G\times G)_{comm})^G$, cf. Section \ref{sec:mapping}. Irreducible characters are a basis for $C((G\times G)_{comm})^G$ and their ${\mathbb Z}$-span gives the Grothendieck ring $K_0(_{D(G)}\mathcal M)$, which is  commutative because $D(G)$ is quasitriangular. The action of $\mathrm{SL}_2(\mathbb Z)$ can be thus naturally interpreted as an action on  $C((G\times G)_{comm})^G$  preserving $K_0(_{D(G)}\mathcal M)$. 
\medskip

  For any $(\mathcal{O},\rho)\in \mathcal{M}(G)$  we denote by $\chi_{(\mathcal{O},\rho)}$ the character of the irreducible $D(G)$-module associated with $(\mathcal{O},\rho)$. The explicit identification of $\chi_{(\mathcal{O},\rho)}$ as an element of $C(G\times G)$ is given by setting
\begin{equation}\label{eq:character}\chi_{(\mathcal{O},\rho)}(h,g)=\chi_{(\mathcal{O},\rho)}(\delta_{(h,g)}), \quad\mbox{ for } h,g\in G.\end{equation}
An identification of ${\mathbb C}[\mathcal M(G)]$ with $C((G\times G)_{comm})^G$ obtained using the language  of equivariant sheaves  appears in  \cite[5.1]{ACR}. The character of $V=\bigoplus_{l\in G}V_l\in\,_{D(G)}\mathcal M$ is mapped to the function whose value on commuting $h$ and $g$ is $\mathrm{Tr}(h|_{V_g})$. Due to the identification \eqref{eq:identify}, the identification \eqref{eq:character}  differs from the one in \cite[5.1]{ACR} by application of $j_1$ as in \eqref{eq:j1j2}.

\medskip

Decomposition of tensor products of irreducible representations reads as 
  \begin{align*}
 \chi_{(\mathcal{O}_1,\rho_1)}\bullet \chi_{(\mathcal{O}_2,\rho_2)}\sum_{(\mathcal{O},\rho)\in\mathcal{M}(G)} N_{\rho_1,\rho_2,\mathcal{O}}^{\mathcal{O}_1,\mathcal{O}_2,\rho}\chi_{(\mathcal{O},\rho)}
\end{align*}
where the $N_{\rho_1,\rho_2,\mathcal{O}}^{\mathcal{O}_1,\mathcal{O}_2,\rho}$ are the so-called {\em fusion coeffiecients}, that is, the non negative integers encoding the multiplicity of the irreducible $D(G)$-module corresponding to $(\mathcal{O},\rho)$ in the tensor product of the irreducible $D(G)$-modules corresponding to  $(\mathcal{O}_1,\rho_1)$ and $(\mathcal{O}_2,\rho_2)$.

\bigskip
 
 By duality, $C((G\times G)_{comm})^G$ is a $k$-subalgebra for $C(G\times G)$ also for the product in $D(G)$. In fact, it is precisely the center of $D(G)$, see \eqref{eq:center}. Orthogonality of characters implies that  
\begin{align*}
  \chi_{(\mathcal{O}_1,\rho_1)}\cdot \chi_{(\mathcal{O}_2,\rho_2)}=\delta(\mathcal{O}_1,\mathcal{O}_2)\delta(\rho_1,\rho_2)\frac{|G|}{|\mathcal{O}_1|\dim(\rho_1)} \chi_{(\mathcal{O}_1,\rho_1)}.
  \end{align*}
 Note that for both  multiplications in $C(G\times G)$ the structure constants of the multiplication of characters are integers, therefore also $K_0(_{D(G)}\mathcal M)$ can be endowed with two different ring structures.
 
 \medskip
 
 The action of $\SFT:=s^{-1}$ on $C((G\times G)_{comm})^G$ satisfies $\SFT^4=1$ and the convolution properties
\begin{align}\label{conv}
    &\SFT(f_1\cdot f_2)=\SFT(f_1)\bullet \SFT(f_2),&\SFT(f_1\bullet f_2)=\SFT(f_1)\cdot \SFT(f_2)
\end{align}
for any $f_1,f_2\in C((G\times G)_{comm})$, motivating the interpretation as a Fourier transform.

Combining \eqref{conv}  with the orthogonality of characters with respect to the multiplication $\cdot$ allows to deduce that the inverse Fourier transform, that is, the $s$-action, diagonalizes the fusion rules. In other words,  $\{\SFT^{-1}\chi_{(\mathcal{O},\rho)}\}_{(\mathcal{O},\rho)\in\mathcal{M}(G)}$ is a basis of simultaneous eigenvectors for the endomorphisms of $C((G\times G)_{comm})^G$ given by $\{\chi_{(\mathcal{O},\rho)}\bullet - \}_{(\mathcal{O},\rho)\in\mathcal{M}(G)}$. Moreover this property allowed to recover in \cite[Theorem 5.3]{tom} the Verlinde formula \cite{V}, an expression of the fusion coefficients in terms of the entries of the matrix of $\SFT$ with respect to the basis of the irreducible characters:
\[N_{\rho_1,\rho_2,\mathcal{O}}^{\mathcal{O}_1,\mathcal{O}_2,\rho}=\sum_{(\mathcal{O}',\rho')\in\mathcal{M}(G)}\frac{\SFT_{(\mathcal{O}',\rho'),(\mathcal{O}_1,\rho_1)}\SFT_{(\mathcal{O}',\rho'),(\mathcal{O}_2,\rho_2)}\overline{\SFT_{(\mathcal{O}',\rho'),(\mathcal{O},\rho)}}}{\SFT_{(\{1\},triv_G),(\mathcal{O}',\rho')}},\]
 see also the account in \cite[8.4.3]{br}.

A proof that the (inverse) Fourier transform diagonalizes the fusion rules for arbitrary finite-dimensional ribbon Hopf algebras was later obtained in \cite{CW}. The transform \eqref{eq:Fourier2} was also related to fusion  data, \cite[Proposition 1.6]{exotic}. The matrix in  \eqref{eq:Fourier} and the diagonal matrix given by $\chi_{\rho}(g)/\chi_\rho(1_G)$ on the pair $(\mathcal O_g,\rho)$ combine to give a representation of $\mathrm{SL}_2(\mathbb Z)$ that actually factors through $\mathrm{PSL}_2(\mathbb Z)$, as the non-abelian Fourier transform matrix is involutive, \cite[1.2]{exotic}. 

\bigskip

In the context of conformal field theory, Verlinde derived a formula  counting the number of group morphisms from fundamental groups of compact Riemann surfaces to $G$, up to $G$-action by conjugation on the target, in terms of characters of centralizers of elements in $G$, \cite{V},\cite[Theorem 4.1]{mason}. If the surface has genus one, it directly implies that $\mathcal M(G)$ is in bijection with the set of $G$-orbits of group homomorphisms $\mathrm{Hom}({\mathbb Z}\times{\mathbb Z},G)^G$. This can also be directly seen, as it amounts to counting pairs of commuting elements in $G$ up to simultaneous conjugation. 
Through the bijection between $\mathrm{Hom}(\mathbb Z\times\mathbb Z,G)^G$ and ${\mathcal M}(G)$ the natural commuting involutions $i_1$ and $i_2$ given by $(\Oc_g,\rho)\mapsto (\Oc_{g^{-1}},\rho)$ and $(\Oc,\rho)\mapsto(\Oc_g,\rho^*)$, \cite[1.5]{exotic} correspond to pre-composition with $(x,y)\mapsto (-x,y)$ and  $(x,y)\mapsto (x,-y)$ on elements in $\mathrm{Hom}({\mathbb Z}\times{\mathbb Z},G)^G$.  
Under the identification $\mathbb C\,\mathcal M(G)\simeq C((G\times G)_{comm})^G$ through characters, the involutions $i_1$ and $i_2$ correspond respectively to the action of the matrices $j_1$ and $j_2$ from \eqref{eq:j1j2}. Indeed an explicit computation of character values yields \cite[3.5]{tom}
\[\chi_{(\mathcal O_l, \rho)}(h,g)=|G| |C_G(l)|\int_{G}\int_{C_G(l)}\delta_{h,g}(h_1lh_1^{-1},h_1g_1h_1^{-1}){\rm Tr}(\rho(g_1))dh_1 \ dg_1. \]
Hence for any $(\mathcal O_l, \rho)\in\mathcal{M}(G)$ there holds
\begin{align*}
    \chi_{(\mathcal O_{l^{-1}}, \rho)}(h,g)=\chi_{(\mathcal O_{l}, \rho)}(h^{-1},g)=j_1\chi_{(\mathcal O_{l}, \rho)}(h,g),
    && \chi_{(\mathcal O_{l}, \rho^*)}(h,g)=\chi_{(\mathcal O_{l}, \rho^*)}(h,g^{-1})=j_2\chi_{(\mathcal O_{l}, \rho)}(h,g).
\end{align*}

The Fourier transform corresponds to the  additional symmetry in $\mathrm{Hom}(\mathbb Z\times\mathbb Z,G)^G$ coming from pre-composing homomorphisms with $(x,y)\mapsto (y,x)$. Indeed, it was observed in \cite[Section 5]{tom} that Lusztig's Fourier transform corresponds to the action of the matrix
 $\begin{pmatrix} 0&1\\1&0\end{pmatrix}\in \mathrm{GL}_2(\mathbb Z)$, that is, to the map $f(h,g)\mapsto f(g,h)$ and that $j_1$ and $j_2$ link Lusztig's Fourier transform with the Fourier transform $\mathcal{S}$. Combined with \eqref{eq:antipode-involution} and our identifications, we obtain the identity
$\mathcal{S}={FT}_G\circ (S_{D(G)^*}\otimes \id)={FT}_G\circ(\id\otimes S_{D(G)})$.

An interpretation of $FT_G$ in terms of symmetry in the variables,  for $G$ a canonical quotient, appears also in \cite[Lemma 5.1]{ACR}, where the authors give credit to \cite{dm90}.

\section{Acknowledgements}
N.C. acknowledges support of ``Ricerca di base - Unipg''.
G.C. and E. C. are members of the INdAM group GNSAGA. The contribution of G. C and E.C. was partially supported by the  
project funded by the European Union – Next Generation EU under the National
Recovery and Resilience Plan (NRRP), Mission 4 Component 2 Investment 1.1 -
Call PRIN 2022 No. 104 of February 2, 2022 of Italian Ministry of
University and Research; Project 2022S8SSW2 (subject area: PE - Physical
Sciences and Engineering) ``Algebraic and geometric aspects of Lie theory''.


\begin{thebibliography}{AFG}
\bibitem{AA}N. Andruskiewitsch, I. Angiono, On Finite dimensional Nichols algebras of diagonal type, Bull. Math. Sci. 7, 353--573, (2017).
\bibitem{ACG5}N. Andruskiewitsch, G. Carnovale, G. A. Garc\'ia. Finite-dimensional pointed Hopf algebras over finite simple groups of Lie type V. Mixed classes in Chevalley and Steinberg groups. Manuscripta Math., 166(3-4) 605–647, (2021).
\bibitem{ACG7}N. Andruskiewitsch, G. Carnovale, G. A. Garc\'ia. Finite-dimensional pointed Hopf algebras over finite simple groups of Lie type VII. Semisimple classes in $\mathrm{PSL}_n(q)$ and $\mathrm{PSp}_{2n}(q)$, J. Algebra, 639 354–397, (2024).
\bibitem{AD}N. Andruskiewitsch, J. Devoto, Extensions of Hopf algebras, Algebra i Analiz 7,  17--52 (1995). 
\bibitem{AFGV1}N. Andruskiewitsch, F. Fantino, M. Gra\~na, L. Vendramin. Finite-dimensional pointed Hopf algebras with alternating groups are trivial. Ann. Mat. Pura Appl. (4), 190(2), 225–245, (2011).
\bibitem{AFGV2}N. Andruskiewitsch, F. Fantino, M. Gra\~na, L. Vendramin. Pointed Hopf algebras over some sporadic simple groups. J. Algebra, 325, 305–320, (2011).
\bibitem{AFS}N. Andruksiewitsch, W. Ferrer Santos, The beginnings of the theory of Hopf algebras, Acta Appl. Math. 108 (2009)  3--17.
\bibitem{AHV}N. Andruksiewitsch, I. Heckenberger, L. Vendramin, Pointed Hopf algebras of odd dimension and Nichols algebras over solvable groups, preprint, arXiv:2411.02304
\bibitem{AN1} N. Andruskiewitsch, S. Natale, Plancherel identity for semisimple Hopf algebras, Commun. Alg. 25 (10), 3239--3254 (1997).  
\bibitem{AN2} N. Andruskiewitsch, S. Natale, Harmonic analysis on semisimple Hopf algebras, Algebra i Analiz 12(5), 3--27, (2000).
\bibitem{AS} N. Andruksiewitsch, H.-J. Schneider,  On the classification of finite-dimensional pointed Hopf algebras, Ann. Math. 171(1), (2010), 375--417.
 \bibitem{ass}I. Assem, D. Simson, A. Skowro\'nski, Elements of the Representation Theory of Associative Algebras, Volume  I, Techniques of Representation Theory, London Mathematical Society Stedent Texts 65, Cambridge University Press, (2006).  
\bibitem{ACR}A.-M. Aubert, D. Ciubotaru, B. Romano, A nonabelian Fourier transform for tempered unipotent representations, Compositio Math., to appear.
\bibitem{Ba}P. Bantay, Orbifolds, Hopf algebras, and the moonshine, Lett. Math. Phys. 22, 187--194 (1991).
\bibitem{B}R. E. Borcherds, Monstrous moonshine and monstrous Lie superalgebras, Invent. Math. 109, 405--444 (1992).
\bibitem{br} M. Brou\'e, On Characters of Finite Groups, Mathematical Lectures from Peking University, Springer, Berlin, (2017).
\bibitem{Br} M. Brou\'e, From Rings and Modules to Hopf Algebras. One Flew Over the Algebraist's Nest, Springer Nature Switzerland, (2024).
\bibitem{camishe}S. Caenepeel, G. Militaru, Zhu Shenglin, Crossed modules and Doi-Hopf modules, Israel Journal of Mathematics 100, 221-247 (1997). 
\bibitem{CC} G. Carnovale, M. Costantini, Finite-dimensional pointed Hopf algebras over finite simple groups of Lie type VI. Suzuki and Ree groups. J. Pure Appl. Algebra, 225(4) Paper No. 106568, (2021).
\bibitem{CP} V. Chari, A. Pressley, A Guide to Quantum Groups, Cambridge University Press, 1994. 
\bibitem{CR} C. Cibils, M. Rosso, Alg\`ebres des chemins quantiques, Adv. Math. 125, 171--199 (1997).
\bibitem{CO} D. Ciubotaru, E. Opdam, On the elliptic nonabelian Fourier transform for unipotent representations of p-Adic groups, In: Representation Theory, Number Theory, and Invariant Theory, Progress in Mathematics, 87--113, Springer International Publishing (2017).
\bibitem{CW}M. Cohen, S. Westreich, Fourier transforms for Hopf algebras, in: Quantum Groups, 115--133,  Contemp. Math., 433
Israel Math. Conf. Proc. American Mathematical Society, Providence, RI, (2007).
\bibitem{CN}J. H. Conway, S.P. Norton,  Monstrous Moonshine, Bull. London Math. Soc., 11(3) 308--339,  (1979).
\bibitem{Coq}R. Coquereaux, Character tables (modular data) for Drinfeld doubles of finite groups,  Proceedings, 7th International Conference on Mathematical Methods in Physics (ICMP2012), Rio de Janeiro, Brazil, April 16-20, 2012,  http://pos.sissa.it/archive/conferences/175/024/ICMP
\bibitem{dancer} K. A. Dancer, P. S. Isac, J. Links,   Representations of the quantum doubles of finite group algebras and spectral parameter dependent solutions of the Yang–Baxter equation, Journal of Mathematical Physics 47, 103511 (2006).
\bibitem{KDC}K. De Commer, Induction for representations of coideal doubles, with an application to quantum ${\mathrm SL}(2,\mathbb R)$, arXiv:2409.18551 (2024).
\bibitem{dm} F. Digne, J. Michel, Fonctions L des variet\'es de Deligne–Lusztig et descente de Shintani, M\'em. S.M.F., vol. 20, Soc. Math. France, Paris, (1985).
\bibitem{dm90}F. Digne, J. Michel, On Lusztig’s parametrization of characters of finite groups of Lie type, Ast\'erisque 181-182, 6,  113--156, (1990).
\bibitem{DM} F. Digne, J. Michel, Representations of Finite Groups of Lie Type, Second Edition, London Mathematical Society Student Texts 95, Cambridge University Press, (2020).
\bibitem{Dijk}R. Dijkgraaf, C. Vafa, E.  Verlinde, H. Verlinde,
The operator algebra of orbifold models, Comm. Math. Phys. 123(3), 485--526, (1989).
\bibitem{DPP} R. Dijkgraaf, V. Pasquier, P. Roche. Quasi Hopf algebras, group cohomology and orbifold models. Nucl. Phys. B (Proc. Suppl.), 18B:60–72, 1990.
\bibitem{ditom} M.S. Dijkhuizen, T.H. Koornwinder, Quantum homogeneous spaces, duality and the quantum 2-spheres, Geom. Ded. 52, 291-315 (1994).
\bibitem{doi}Y. Doi, Unifying Hopf modules. J. Algebra 153 (1992), 373–385.
\bibitem{drinfeld} V.G. Drinfel’d, Quantum groups, in: Proceedings of the I.C.M., Berkeley, (1986), American
Math. Soc., 1987, 798–820.
\bibitem{etingof}P. Etingof, On some properties of quantum doubles of finite groups, J. of Algebra 394, 1--6 (2013).
\bibitem{egno} P. Etingof, S. Gelaki, V. Ostrik, D. Nykshych, Tensor Categories, Mathematical Surveys and Monographs 205, AMS, 2015.
\bibitem{erw}P. Etingof, E. C. Rowell, S. J. Witherspoon, Braid group representations from quantum doubles of finite groups, Pacific J. Math. 234(1), 33--41, (2008).
\bibitem{FF} J. Fjelstad, J. Fuchs,  Mapping class group representations from Drinfeld doubles of finite groups, Journal Knot theory and Ramifications, 29(05) 2050033, (2020).
\bibitem{gould}M. Gould, Quantum double finite group algebras and their representations,  Bulletin of the Australian Mathematical Society 48(02), 275--301, (1993).
\bibitem{tohoku}A. Grothendieck, Sur quelques points d'alg\`ebre homologique, T\^{o}hoku Mathematical Journal (2), 9(2), 119--221, (1957).
\bibitem{gust}P. P.  Gustafson, Finiteness for mapping class group representations from twisted Dijkgraaf–Witten theory,  Journal of Knot Theory and Its Ramifications 27(6), 1850043 (2018). 
\bibitem{HS}I. Heckenberger, H.-J. Schneider, Hopf algebras and Root systems,  Mathematical Surveys and Monographs 247, American Mathematical Society (2020). 
\bibitem{lam-rad} L. A. Lambe, D. Radford, Algebraic aspects of the quantum Yang-Baxter equation, J. Algbra 154, 228-288 (1993).
\bibitem{js} A. Joyal, R. Street, Tortile Yang-Baxter operators in tensor categories. Journal of Pure and Applied Algebra, 71(1), 43–51, (1991).
\bibitem{JS}A. Joyal, R. Street, Braided Tensor Categories, Adv. Math. 102, 20--78, (1993).
\bibitem{Kac} G. I. Kac, A generalization of the principle of duality for groups. Dokl. Akad. Nauk SSSR, 138,  275–278, (1961).
\bibitem{KP}G. I. Kac, V. G. Paljutkin, Finite ring groups, Trudy Moskov. Mat. Ob\v{s}\v{c}. 15 (1966), 224–261.
\bibitem{Ke}M. Keilberg, Quasitriangular structures of the double of a finite group, Comm. Algebra 46, 5146-5178 (2018).
\bibitem{ke2}M. Keilberg, Automorphisms of the Doubles of Purely Non-abelian Finite Groups, Algebr Represent Theory 18, 1267--1297, (2015).
\bibitem{kassel}C. Kassel, Quantum Groups, GTM 155, Springer-Verlag New-York, 1995.
\bibitem{kitaev} A. Yu. Kitaev, Fault-tolerant quantum computation by anyons, Annals of Physics 303, 2--30, (2003).
\bibitem{tom-muller} T.H. Koornwinder and N.M. Muller. The quantum double of a (locally) compact group. J. Lie Theory, 7:33–52, (1997) and 8:187, (1998) (erratum).
\bibitem{tom} T.H. Koornwinder, B. J. Schroers, J. K. Slingerland, F. Bais, Fourier transform and the Verlinde formula for the quantum double of a finite group, Journal of Physics A Mathematical and General 32(48),8539--8549, (1999).
\bibitem{tom-tensor}T. H. Koornwinder, F. Bais, N.M. Muller, Tensor product representations of the quantum double of a compact group, Comm. Math. Phys. 198(1) 157--186, (1998).
\bibitem{kop} M. Koppinen, Uncertainty inequalities for Hopf algebras, Comm. Algebra 22(3), 1083--1101, (1994).
\bibitem{koppinen} M. Koppinen, Variations on the smash product with applications to group-graded rings. J. Pure Appl. Algebra 104, 61–80  (1995).

\bibitem{Orange} G. Lusztig, Characters of Reductive Groups over a Finite Field, Princeton University Press (1984).
\bibitem{arcata} G. Lusztig, Leading coefficients of character values of Hecke algebras,in the Arcata Conference on Representations of Finite Groups, P.Fong ed., Proceedings of Symposia in Pure Mathematics, vo147 (A.M.S, 1987).
\bibitem{Lu-QG}G. Lusztig, Introduction to Quantum Groups, Birkh\"auser, (1993).
\bibitem{exotic}G. Lusztig, Exotic Fourier Transform, Duke Mathematical Journal 73(1),  227--241, (1994).
\bibitem{ly} V. V. Lyubashenko, Ribbon abelian categories as modular categories, J. Knot Theory and its Ramific. 5, 311--403, (1996).
\bibitem{LM}V. V. Lyubashenko, S. Majid, Braided groups and quantum Fourier transform, J. Algebra 166(3), 506--528, (1994).
\bibitem{Maj}S. Majid, Doubles of quasitriangular Hopf algebras, Comm. Algebra 19(11), 3061-3073 (1991).
\bibitem{maj-center}S. Majid, Representations, duals and quantum doubles of monoidal categories, Rend. Circ. Mat. Palermo (2) 26, 197--206, (1991). 
\bibitem{mabook}S. Majid, Foundations of Quantum Group Theory, Cambridge University Press, (1995). 
\bibitem{maj-qm} S. Majid, The quantum double as quantum mechanics. J. Geom. Phys., 13, 169--202  (1994).
\bibitem{maj-LNPAM} S. Majid, Classification of differentials on quantum doubles and finite non commutative geometry, Lect. Notes Pure Appl. Math. 239, 167-188, Marcel Dekker (2004).
\bibitem{majmc} S. Majid, L. S. McCormack, Quantum geometric Wigner construction for D(G) and braided racks, arXiv:2407.11835.
\bibitem{mason}G. Mason, The quantum double of a finite group and its role in conformal field theory. In: Groups '93 Galway/St. Andrews, 2, 405--417.
London Math. Soc. Lecture Note Ser., 212, Cambridge University Press, Cambridge, (1995).
\bibitem{Matsu}H. Matsumoto, G\'en\'erateurs et relations des groupes de Weyl g\'en\'eralis\'es, C. R. Acad. Sci. Paris, 258, 3419--3422, (1964).
\bibitem{muger}M. M\"uger, Quantum double actions on operator algebras and orbifold quantum field theories,  Comm. Math. Phys. 191(1), 137--181, (1998).
\bibitem{Mu}E. M\"uller, Some topics on Frobenius-Lusztig kernels, I, J. Algebra 206, 624--658 (1998).
\bibitem{nichols}W. D. Nichols, Bialgebras of type one, Comm. Alg. 6,  1521--1552 (1978).
\bibitem{radford} D. Radford, Minimal quasitriangular Hopf algebras, J. of Algebra 157, 285--315 (1993). 
\bibitem{radbook}D. Radford, Hopf Algebras, Series on Knots and Everything Vol. 49, World Scientific, (2012). 
\bibitem{rad-tow}D. Radford, J. Towber, Yetter-Drinfel’d categories
associated to an arbitrary bialgebra, Journal of Pure and Applied Algebra 87, 259--279 (1993).
\bibitem{reshe}N. Y. Reshetikhin, M. A. Semenov-Tian-Shansky, Quantum R-matrices and factorization problems, Journal Geometry Physics 5, 533--550 (1988). 
\bibitem{RT}N. Y. Reshetikhin, V. Turaev, Ribbon graphs and their invariants derived from quantum groups, Comm. Math. Phys 127, 1--26, (1990).
\bibitem{rosso} M. Rosso,  Groupes quantiques et algebres de battage quantiques, C.R.A.S. (Paris) 320,  145--148 (1995).
\bibitem{rosso-inv}M. Rosso, Quantum groups and quantum shuffles, Inventiones Math. 133 399--416, (1998).
\bibitem{rowell2}E. Rowell, Two paradigms for topological quantum computation, Contemp. Math. 482, 165--178, AMS, Providence, RI, (2009).
\bibitem{rowell}E. Rowell, Zhenghan Wang, Mathematics of topological quantum computing, Bull. Amer. Math. Soc.  (New Series) 55(2), 183--238, (2018).
\bibitem{Sch} P. Schauenburg, Hopf modules and Yetter-Drinfeld modules, J. Algebra 169, 874--890, (1994).
\bibitem{Sch2}P. Schauenburg, A characterization of the Borel-like subalgebras of quantum enveloping algebras, Comm. Algebra 24, 2811--2823, (1996).
\bibitem{Segal} G. Segal, The definition of conformal field theory, in: Differential geometrical methods in theoretical physics (Como, 1987), NATO Adv. Sci. Inst. Ser. C Math. Phys. Sci., 250, Kluwer Acad. Publ., Dordrecht, 1988, 165-171.
\bibitem{SS} S. Shnider, S. Sternberg, Quantum Groups - from coalgebras to Drinfeld algebras - a guided tour, Graduate Texts in Mathematical Physics, Volume II,  International Press Publications, Hong Kong, (1993).
\bibitem{Shoji} T. Shoji, Character sheaves and almost characters of reductive groups, I, II. Adv.
Math. 111, 244–313, 314–354, (1995).
\bibitem{springer}T.A. Springer, Sheaf Theory, Lecture Notes, Mathematical Institute, Universty Utrecht, (1977). 
\bibitem{sweedler}M. Sweedler, Hopf Algebras, W. A. Benjiamin, Inc., New York, (1969). 
\bibitem{tu}V.G. Turaev, Quantum Invariants of Knots and 3-Manifolds, Studies in Mathematics 18, de Gruyter, New York (1994).
\bibitem{Va}L. Va\u{\i}nerman, Gel'fand pairs of quantum groups, hypergroups and q-special functions, Applications of Hypergroups and Related Measure Algebras, (Seattle, WA, 1993), Contemp. Math. vol. 183, Amer. Math. Soc. Providence, RI, 373--394, (1995).
\bibitem{V}E. Verlinde, Fusion rules and modular transformations in 2D conformal field theory, Nucl. Phys. B, 300, 360--376 (1988). 
\bibitem{wither}S. Witherspoon, 
The Representation Ring of the Quantum Double of a Finite Group, Journal of Algebra 179, 305-329 (1996).
\bibitem{30}S. L. Woronowicz, Twisted SU (2)-group, an example of a non-commutative differential calculus. Publ. RIMS (Kyoto) 23, 117--181 (1987).
\bibitem{W}S. L. Woronowicz, Compact Matrix Pseudogroups,  Commun. Math. Phys. 111, 613--665 (1987).
\bibitem{wor-diff-cal}S. L. Woronowicz, Differential Calculus
on Compact Matrix Pseudogroups (Quantum Groups), Commun. Math. Phys. 122, 125--170 (1989).
\bibitem{yetter}D. Yetter, Quantum groups and representations of monoidal categories, Math. Proc. Camb. Phil. Soc.  108, 261--290  (1990).
\end{thebibliography}
\end{document}